\documentclass{gtart_h}

\def\ifplaintex{\expandafter\ifx\csname documentclass\endcsname\relax}

\def\gtp{{\mathsurround=0pt\it $\cal G\mskip-2mu$eometry \&\ 
$\cal T\!\!$opology $\cal P\!$ublications}}  

\def\recd{{\small Received:\qua\receiveddate\ifx\reviseddate\relax
\else\qquad Revised:\qua\reviseddate\fi\par}} 

\def\lognumber#1{\def\thelognumber{#1}}
\def\volumenumber#1{\def\thevolumenumber{#1}}
\def\volumeyear#1{\def\thevolumeyear{#1}}
\def\papernumber#1{\def\thepapernumber{#1}}
\def\pagenumbers#1#2{\def\startpage{#1}\def\finishpage{#2}}
\def\published#1{\def\publishdate{#1}}

\def\received#1{\def\receiveddate{#1}}

\def\accepted#1{\def\accepteddate{#1}}
\def\asciititle#1{\def\theasciititle{#1}}
\def\covertitle#1{\def\thecovertitle{#1}}

\long\def\asciiabstract#1{\long\def\theasciiabstract{#1}}
\def\asciikeywords#1{\def\theasciikeywords{#1}}

\let\\\par\let\thelognumber\relax\let\thevolumenumber\relax
\let\thepapernumber\relax\let\thevolumeyear\relax\let\startpage\relax
\let\finishpage\relax\let\publishdate\relax\let\receiveddate\relax
\let\reviseddate\relax\let\accepteddate\relax\let\theasciititle\relax
\let\thecovertitle\relax\let\theasciiauthors\relax
\let\theasciiabstract\relax\let\theasciikeywords\relax

\let\theasciiemail\relax

\ifplaintex
\font\logobig=cmssbx10 scaled 3836
\font\logomed=cmssbx10 scaled 2557
\else
\font\logobig=cmssbx10 scaled 4200
\font\logomed=cmssbx10 scaled 2800
\fi

\long\def\makeagttitle{   
\count0=\startpage
\agt\hfill      
\hbox to 45truept{\vbox to 0pt{\vglue -13truept{\logomed A\kern -.37em{\logobig 
T}\kern -.38em G}\vss}\hss}
\break
{\small Volume \thevolumenumber\ (\thevolumeyear)
\startpage--\finishpage\nl
Published: \publishdate}

\vglue .25truein

{\parskip=0pt\leftskip 0pt plus
1fil\def\\{\par\smallskip}{\Large\bf\thetitle}\par\medskip} \vglue
0.05truein

{\parskip=0pt\leftskip 0pt plus 1fil\def\\{\par}{\sc\theauthors}
\par\medskip}%
 
\vglue 0.03truein 

{\small\leftskip 25truept\rightskip 25truept{\bf Abstract}\stdspace\theabstract

{\bf AMS Classification}\stdspace\theprimaryclass
\ifx\thesecondaryclass\relax\else; \thesecondaryclass\fi\par
{\bf Keywords}\stdspace \thekeywords\par}\vglue 7truept

}   

\ifplaintex
\font\phead=cmsl9 scaled 950
\font\pnum=cmbx10 scaled 913
\font\pfoot=cmsl9 scaled 950
\headline{\vbox to 0pt{\vskip -4.5mm\line{\small\phead\ifnum
\count0=\startpage ISSN 1472-2739 (on-line) 1472-2747 (printed)
\hfill {\pnum\folio}\else\ifodd\count0\def\\{ }%
\ifx\theshorttitle\relax\thetitle\else\theshorttitle\fi\hfill{\pnum\folio}
\else\def\\{ and }{\pnum\folio}\hfill\ifx\theshortauthors\relax\theauthors
\else\theshortauthors\fi\fi\fi}\vss}}
\footline{\vbox to 0pt{\vglue 0mm\line{\small\pfoot\ifnum\count0=\startpage
\copyright\ \gtp\hfill\else
\agt, Volume \thevolumenumber\ (\thevolumeyear)\hfill\fi}\vss}}
\else
\font=cmsl9 scaled 1050
\font\lnum=cmbx10 
\font=cmsl9 scaled 1050
\makeatletter
\def\@oddhead{{\small\lhead\ifnum\count0=\startpage ISSN 1472-2739 
(on-line) 1472-2747 (printed)\hfill {\lnum\number\count0}\else\ifodd\count0
\def\\{ }\ifx\theshorttitle\relax \thetitle \else\theshorttitle\fi\hfill
{\lnum\number\count0}\else\def\\{ and }{\lnum\number\count0}
\hfill\ifx\theshortauthors\relax 
\theauthors\else\theshortauthors\fi\fi\fi}}\def\@evenhead{\@oddhead}
\def\@oddfoot{\small\lfoot\ifnum\count0=\startpage\copyright\ \gtp\hfill\else
\agt, Volume \thevolumenumber\ (\thevolumeyear)\hfill\fi}
\def\@evenfoot{\@oddfoot}
\makeatother
\fi
\let\maketitlepage\makeagttitle
\let\makeshorttitle\maketitlepage
\let\maketitle\maketitlepage

\newwrite\gtoutfile
\long\gdef\makeheadfile{  
{\def\\{, }\def\s{ }
\immediate\openout\gtoutfile head.xxx
\immediate\write\gtoutfile{Proxy-for: \ifx\theasciiauthors\relax
\theauthors\else\theasciiauthors\fi\s<\ifx\theasciiemail\relax\theemail\else\theasciiemail\fi>}
\immediate\write\gtoutfile{\noexpand\\}
\immediate\write\gtoutfile{Authors: \ifx\theasciiauthors\relax
\theauthors\else\theasciiauthors\fi}
{\def\\{ }\immediate\write\gtoutfile{Title: \ifx\theasciititle\relax
\thetitle\else\theasciititle\fi}}
\immediate\write\gtoutfile{Subj-class: GT or SG, GR etc}
\immediate\write\gtoutfile{MSC-class: \theprimaryclass\ifx\thesecondaryclass\relax\else, \thesecondaryclass\fi}
\immediate\write\gtoutfile{Journal-ref: Algebr. Geom. Topol. \thevolumenumber\s
(\thevolumeyear) \startpage-\finishpage}
\immediate\write\gtoutfile{Comments: Published by Algebraic and
Geometric Topology at}
\immediate\write\gtoutfile{\s\s\s  http://www.maths.warwick.ac.uk/agt/AGTVol\thevolumenumber/agt-\thevolumenumber-\thepapernumber.abs.html}
\immediate\write\gtoutfile{\noexpand\\}
\immediate\write\gtoutfile{}
\ifx\theasciiabstract\relax
\immediate\write\gtoutfile{\theabstract}\else
\immediate\write\gtoutfile{\theasciiabstract}\fi
\immediate\write\gtoutfile{}
\immediate\write\gtoutfile{\noexpand\\}
\immediate\write\gtoutfile{}
\immediate\closeout\gtoutfile}}  

\def\maketitlepage{\makeagttitle\makeheadfile}
\let\makeshorttitle\maketitlepage
\let\maketitle\maketitlepage

\lognumber{30}
\volumenumber{4}
\volumeyear{2004}
\papernumber{30}
\pagenumbers{647}{684}
\received{3 December 2002} 
\accepted{9 August 2004}
\published{3 September 2004}

\usepackage{amsmath, amssymb, amscd}

\newtheorem{thm}{Theorem}[section]
\newtheorem{lm}[thm]{Lemma}
\newtheorem{cor}[thm]{Corollary}
\newtheorem{pro}[thm]{Proposition}
 
\theoremstyle{definition}
\newtheorem{df}[thm]{Definition}

\numberwithin{equation}{section}

\let\Bbb\mathbb
\def \BR {\Bbb{R}}
\def \BC {\Bbb{C}}
\def \BZ {\Bbb{Z}}

\def \J {\cal J}
\def \P {\cal P}
\def \CM {\cal M}
\def \hM {\hat{{\cal M}}}

\def \la {\langle}
\def \ra {\rangle}
\def \a {\alpha}

\def \co {\colon\thinspace}
\def \g {\gamma}

\def \d {\delta}

\def \Lam {\Lambda}
\def \lam {\lambda}

\def \n {\nabla}
\def \o {\omega}
\def \op {\oplus}
\def \ov {\overline}
\def \p {\phi}
\def \r {\rho}
\def \s {\sigma}
\def \Si {\Sigma}
\def \t {\tau}
\def \ti {\tilde}
\def \T {\Theta}

\def \O {\Omega}

\def \bd {\partial}
\def \pj {\overline{\partial}_J}

\def \x {\times}
\def \ox {\otimes}
\def \ve {\varepsilon}
\def \D {\mathbf{D}}

\begin{document}

\title[Lagrangian embeddings]
{The ${\BZ}$--graded symplectic Floer cohomology\\of monotone 
Lagrangian sub--manifolds }
\covertitle{The ${\noexpand\bf Z}$--graded symplectic Floer cohomology\\of monotone 
Lagrangian sub--manifolds }
\asciititle{The Z-graded symplectic Floer cohomology of monotone 
Lagrangian sub-manifolds }

\author{Weiping Li}
\address{Department of Mathematics, Oklahoma State University\\Stillwater, 
Oklahoma 74078-0613}
\email{wli@math.okstate.edu}

\begin{abstract} We define an integer graded symplectic Floer cohomology
and a Fintushel--Stern type spectral sequence which are new invariants
for monotone Lagrangian sub--manifolds and exact isotopes. The 
${\BZ}$--graded symplectic Floer 
cohomology is an integral lifting of the usual ${\BZ}_{\Si (L)}$--graded 
Floer--Oh
cohomology. We prove the K\"{u}nneth formula for
the spectral sequence and an ring structure on it. The ring structure on the
${\BZ}_{\Si (L)}$--graded Floer cohomology is induced from the ring
structure of the cohomology of the Lagrangian sub--manifold via the spectral 
sequence.
Using the ${\BZ}$--graded symplectic Floer
cohomology, we show some intertwining relations among the Hofer energy
$e_H(L)$ of the embedded Lagrangian, the minimal symplectic action $\s(L)$, the
minimal Maslov index $\Si (L)$ and the smallest integer $k(L, \p)$ of the
converging spectral sequence of the Lagrangian $L$.
\end{abstract}

\asciiabstract{We define an integer graded symplectic Floer cohomology
and a Fintushel-Stern type spectral sequence which are new invariants
for monotone Lagrangian sub-manifolds and exact isotopes. The Z-graded
symplectic Floer cohomology is an integral lifting of the usual
Z_Sigma(L)-graded Floer-Oh cohomology. We prove the Kunneth formula
for the spectral sequence and an ring structure on it. The ring
structure on the Z_Sigma(L)-graded Floer cohomology is induced from
the ring structure of the cohomology of the Lagrangian sub-manifold
via the spectral sequence.  Using the Z-graded symplectic Floer
cohomology, we show some intertwining relations among the Hofer energy
e_H(L) of the embedded Lagrangian, the minimal symplectic action
sigma(L), the minimal Maslov index Sigma(L)$ and the smallest integer
k(L, phi$ of the converging spectral sequence of the Lagrangian L.}

\asciikeywords{Monotone Lagrangian sub-manifold,
Maslov index, Floer cohomology, spectral
sequence}

\maketitle

\section{Introduction}

In this paper, we construct
a ${\BZ}$--graded symplectic Floer cohomology of monotone Lagrangian 
sub-manifolds by a completely algebraic topology method. 
This is a local symplectic invariant in terms of symplectic
diffeomorphisms. We show that there
exists a spectral sequence which converges to a global symplectic invariant
(the ${\BZ}_{\Si (L)}$--graded Floer cohomology, where 
${\BZ}_{\Si (L)}$ is the 
minimal Maslov number of the monotone Lagrangian sub--manifold $L$).
The ${\BZ}$--graded symplectic Floer cohomology is an integral
lifting of the ${\BZ}_{\Si (L)}$--graded symplectic Floer cohomology.
By exploiting the properties of our ${\BZ}$--graded 
symplectic Floer cohomology, 
we show that there is a relation between the ${\BZ}$--graded symplectic
Floer cohomology and 
the restricted symplectic
Floer cohomology constructed in \cite{ch} (see \S 5) via Hofer's
symplectic energy. This may give an interesting way to understand
the Hofer symplectic energy through the ${\BZ}$--graded symplectic
Floer cohomology.
We borrow some ideas from the instanton Floer 
theory, and construct the Fintushel-Stern type spectral sequence in
\cite{fs2} for
monotone Lagrangian sub--manifolds. Our 
method is in the nature of algebraic topology (see \cite{sp}), and is very
different from the method in \cite{oh2} which is used local Darboux
neighborhoods.
We hope that there will be more algebraic cohomology operations which
can be induced to the symplectic Floer cohomology through the quantum effects
of higher differentials in our spectral sequence.

Let $(P, \o)$ be a monotone symplectic manifold and $L$ 
be a monotone Lagrangian
sub-manifold in $(P, \o)$. Let $Z_{\p}$ be the critical point set of the
symplectic action $a_{\p}$ (see \S 2), where $\p \in \text{Symp}_0(P)$
is a symplectic diffeomorphism generated by a time--depended 
Hamiltonian function. The set $\text{Im}\, (a_{\p})(Z_{\p})$ is discrete.
For $r \in {\BR } \setminus \text{Im}\, (a_{\p})(Z_{\p})
= {\BR }_{L, \p}$, we can associate
the ${\BZ}_2$--modules $I^{(r)}_*(L, \p; P)$ with an
integer grading. 
The ${\BZ}$--graded symplectic Floer cohomology $I^{(r)}_*(L, \p; P)$ 
depends on $r$:

(i)\qua if $[r_0, r_1] \subset {\BR }_{L, \p}$, then
$I^{(r_0)}_*(L, \p; P) = I^{(r_1)}_*(L, \p; P)$; 

(ii)\qua $I^{(r)}_{* + \Si (L)}
(L, \p; P) = I^{(r+ \s (L))}_*(L, \p; P)$, where 
$\Si (L) ( > 0)$ is the minimal Maslov number of $L$ and
$\s (L) ( > 0)$ is the minimal 
number in $\text{Im}\, I_{\o }|_{\pi_2 (P, L)}$ (see Definition~\ref{d2.2}
for $\Si (L)$ and $\s (L)$).

Our main results are the following theorems.
\medskip

\noindent{\bf Theorem A}\qua {\sl Let $L$ be a  
monotone Lagrangian sub--manifold in $(P, \o)$. If $\Si (L) \geq 3$, then

{\rm(1)}\qua there exists an isomorphism
\[\p_{01}^n \co I^{(r)}_n(L, \p^0; P, J^0) \to I^{(r)}_n(L, \p^1; P, J^1), \]
for $n \in {\BZ}$, and a continuation 
$(J^{\lam }, {\p }^{\lam })_{0\leq \lam \leq 1} \in {\P }_1$ which
is regular at the ends.

{\rm(2)}\qua there is a spectral sequence $(E^k_{n, j}, d^k)$ with
\[ E^1_{n,j}(L, \p; P, J) \cong I_n^{(r)}(L, \p; P, J), \ \ \
n\equiv j \ \ (mod \ \ \Si (L)), \]
\[ d^k \co E^k_{n, j}(L, \p; P, J) \to
E^k_{n + \Si (L) k +1, j +1}(L, \p; P, J) , \]
and
\[E^{\infty }_{n,j}(L, \p; P, J) \cong F^{(r)}_n HF^j(L, \p; P, J)/
F^{(r)}_{n + \Si (L)}HF^j(L, \p; P, J) . \]
}
\medskip

\noindent{\bf Theorem B}\qua 
{\sl For $\Si (L) \geq 3$ and $k \geq 1$, 
$E^k_{n, j}(L, \p; P, J)$ is a symplectic invariant
under continuous deformations of $(J^{\lam }, {\p }^{\lam })$ within the
set of continuations.}
\medskip

For $k \geq 1$ and $r \in {\BR }_{L, \p}$, all the $E^k_{*, *}(L, \p; P, J)$
are new symplectic invariants. They provide potentially interesting 
invariants for the symplectic topology of $L$. 
Let $k(L, \p)$ be the minimal $k$ for
which $E^k_{*, *}(L, \p) = E^{\infty }_{*, *}(L, \p)$. 
So $k(L, \p)$ is an invariant of $(L, \p)$.
Using the ${\BZ}$--graded symplectic Floer cohomology and 
the spectral sequence in Theorem A,
we obtain the K\"{u}nneth formulae for each term of the spectral sequence
with ${\BZ}_2$--coefficients. 
Since we work on the ${\BZ}_2$--coefficients, there are K\"{u}nneth formulae
for the induced spectral sequence of the product monotone Lagrangian  
sub-manifold $(L_1 \x L_2, \p_1 \x \p_2)$.
Then we study the Poincar\'{e}--Laurent polynomial for $L_1 \x L_2$
in terms of the Poincar\'{e}--Laurent polynomials for $(L_i, \p_i) (i=1, 2)$.
For certain Lagrangian imbeddings, we obtain an internal cup--product structure
on the spectral sequence which is descended from the usual cup product
of the cohomology $H^*(L, {\BZ}_2)$, by studying the 
$H^*(L, {\BZ}_2)$--module
structure on the spectral sequence and the 
${\BZ}_2$--coefficients. The index
in the internal product (\ref{uuu})
is unusual from the internal product structure due to the Maslov index shift.
From the quantum effect aroused from the higher differentials on
$I_*^{(r)}(L, \p; P, J)$ in the
spectral sequence in Theorem A, the ring
$(HF^{*-m}(L, \p; P), \cup_{\infty})$ on the ${\BZ}_{\Si (L)}$--graded
symplectic Floer cohomology can be thought of as
the quantum effect of the cohomology ring
$(H^*(L; {\BZ}_2), \cup)$ (see \S 5.2).
Note that our cup-product structure  is different with the multiplicative
structure defined in \cite{fooo, lt, rt}. In \cite{lt, rt}, 
the Floer cohomology is the cohomology of the symplectic manifold, 
only the cup-product structure is deformed, 
i.e., the same cohomology group with different ring structures is
studied in \cite{lt, rt}. Our induced cup--product on $E_{*, *}^k(L, \p; P, J)$
may well have that the cohomology groups 
are different from the cohomology of the Lagrangian sub--manifolds
(see \cite{li} for instance).
\medskip

\noindent{\bf Theorem C}\qua {\sl {\rm(1)}\qua
For the monotone Lagrangian $L_1 \x L_2$ in $(P_1 \x P_2, \o_1 \op \o_2)$ with
$I_{\o_i} = \lam I_{\mu, L_i}$ and $\Si (L_i) = \Si (L) \geq 3$ ($i=1, 2$),
we have, for $k \geq 1$,
$$E^k_{n,j}(L_1 \x L_2, \p_1 \x \p_2; P_1 \x P_2) \cong$$ 
\begin{equation}
\bigoplus_{ n_1 + n_2 = n, j_1 + j_2 = j \pmod {\Si (L)} }
E^k_{n_1,j_1}(L_1, \p_1; P_1) \ox E^k_{n_2,j_2}(L_2, \p_2; P_2).
\end{equation}
{\rm(2)}\qua For the monotone Lagrangian $L \stackrel{i}{\hookrightarrow} P$
with $i^* \co 
H^*(P; {\BZ}_2) \to H^*(L; {\BZ}_2)$ surjective and $\Si (L) \geq 3$,
the spectral sequence $E^k_{*-m, *}(L, \p; P)$ carries an ring structure
which is descended from the cohomology ring $(H^*(L; {\BZ}_2), \cup )$.}

As an easy consequence of Theorem C, we obtain a generalization of Theorem
1 and Theorem 3 of \cite{fl1} as stated in Corollary~\ref{inh}. This proves
the Arnold conjecture that the monotone Lagrangian intersections is
bounded below by the ${\BZ}_2$--cuplength of the Lagrangian sub-manifold
(see \S 5.2).
Using a result of Gromov and the Poincar\'{e}--Laurent polynomial associated
to the spectral sequence, we 
show that the
four invariants $\s (L), \Si (L), e_H(L)$ and $k(L, \p)$ play
important roles in the ${\BZ}$--graded symplectic Floer cohomology and
the study of Lagrangian embeddings in \S 5.1. 
We obtain Chekanov's result by using the ${\BZ}$--graded 
symplectic Floer cohomology. 
Our study suggests a possible relation between
the ${\BZ}$--graded 
symplectic Floer cohomology and Hofer's symplectic energy for monotone
Lagrangian sub-manifolds. In fact we conjecture
that Hofer's symplectic energy of a monotone Lagrangian sub--manifold $L$
with $\Si (L) \geq 3$ is a positive multiple of $\s (L)$ 
(More precisely, $e_H(L) = (k(L, \p)-1) \s (L)$). 
We will discuss this problem
elsewhere. It would be also
interesting to link the ${\BZ}$--graded symplectic Floer cohomology with
the (modified) Floer cohomology with 
Novikov ring coefficients in \cite{fooo}. 

The paper is organized as follows. In \S 2, we define the ${\BZ}$--graded 
symplectic Floer cohomology for monotone Lagrangian 
sub--manifolds. 
Its invariance under the symplectic continuations  
is given in \S 3. Theorem A (1)
is proved in \S 3. Theorem A (2), Theorem B and Theorem C (1)
(Theorem~\ref{kus}) 
are proved in \S 4. In \S
5, we give some applications related to Chekanov's construction and
Lagrangian embeddings; 
at the last subsection \S 5.2, the proof of Theorem C (2) (Theorem~\ref{cupp})
is given.

\section{The $\BZ$--graded Floer cohomology for Lagrangian
intersections}

In this section, we define the $\BZ$--graded symplectic Floer cohomology,
and discuss some basic properties.

Let $(P, \o)$ be an oriented, connected and compact (or tamed) symplectic 
manifold with a closed non--degenerate
2-form $\o $. The 2--form $\o$ defines the cohomology class 
$[\o ] \in H^2(P, {\BR })$. By
choosing an almost complex structure $J$ on $(P, \o)$ such that
$\o ( \cdot , J \cdot )$ defines a Riemannian metric, we  have an
integer valued cohomology class $c_1(P) \in H^2(P, {\BZ})$ (the first Chern
class). These two cohomology classes define two homomorphisms
\[I_{\o }\co \pi_2(P) \to {\BR }; \ \ \ \ \ \
I_{c_1}\co \pi_2(P) \to {\BZ }, \]
by $I_{\o }(u) = \int_{S^2} u^*(\o)$ and $I_{c_1}(u) = \int_{S^2} u^*(c_1)$
for $u \in \pi_2(P)$.
If $u \colon (D^2, \partial D^2) \to (P, L)$ is a smooth map of pairs, 
up to homotopy, there is 
a unique trivialization of the pull-back bundle $u^* TP
\cong D^2 \x {\BC}^m$. The trivialization of the symplectic vector bundle
defines a map from $S^1 = \bd D^2$ to the set ${\Lam }({\BC}^m)$ of
Lagrangians in ${\BC}^m$. Let $\mu \in H^1({\Lam }({\BC}^m), \BZ)$ 
be the well--known Maslov class. Then the map
$I_{\mu, L}\co \pi_2(P, L) \to {\BZ}$ is defined 
by $I_{\mu, L}(u ) = \mu (\bd D^2)$. The Maslov index is invariant
under any symplectic isotopy of $P$.
\begin{df} (i)\qua The symplectic manifold $(P, \o)$ is {\it monotone} if
$I_{\o} = \a I_{c_1}$ for some $\a \geq 0.$

(ii)\qua A Lagrangian sub--manifold $L$ in $(P, \o)$ is {\it monotone} if
$I_{\o} = \lam I_{\mu, L}$ for some $\lam \geq 0.$
\end{df}
For $\a = 0$ and $\lam = 0$ case, the manifold
$(P, \o)$ and $L$ are monotone defined by
Floer \cite{fl2, fl3}. The notion of monotone Lagrangian
sub-manifolds is introduced by Oh \cite{oh}.
The monotonicity is preserved under the exact deformations of $L$.
By the canonical homomorphism $f\co \pi_2(P) \to \pi_2(P,L)$,
\[ I_{\o}(x) = I_{\o}(f(x)), \ \ \ \ \
I_{\mu , L}(f(x)) = 2 I_{c_1}(x), \]
for $x \neq 0 \in \pi_2(P)$. If the Lagrangian sub-manifold $L$ is 
monotone, 
then $(P, \o)$ is also a monotone symplectic manifold with 
$2 \lam =  \a$. In fact, the constant $\lam$ does not depend
on the Lagrangian $L$, but depends only on the $(P, \o)$ if $I_{\o}|_{
\pi_2(P)} \neq 0$. 

\begin{df}\label{d2.2} (i)\qua 
Define $\s (L)$ to be the positive minimal number in 
the set $\text{Im}\, I_{\o}|_{\pi_2(P, L)} 
\subset \BR$.
Define $\Si (L)$ to be the positive generator for the subgroup
$[\mu|_{\pi_2(P, L)}] = \text{Im}\, I_{\mu, L}$ in $\BZ$. 

(ii)\qua A Lagrangian sub--manifold
$L$ is called rational if $\text{Im}\, I_{\o}|_{\pi_2(P, L)} = \s (L) \BZ$
is a discrete subgroup of $\BR$ and $\s (L) > 0$. 
For a monotone Lagrangian, we
have $\s (L) = \lam \Si (L)$ for some $\lam \geq 0$.
\end{df}
 
Let $H\co P \x {\BR} \to {\BR}$ be a 
smooth real valued function and let $X_H$ be 
defined by $\o (X_H, \cdot ) = dH$. Then the ordinary differential 
equation
\begin{equation} \label{hami} 
\frac{dx}{dt} = X_H(x(t)) ,
\end{equation} 
is called a Hamiltonian equation associated with the time--dependent
Hamiltonian
function $H$, or with the Hamiltonian vector field $X_H$. 
Equation (\ref{hami}) defines a family $\p_{H,t}$
of diffeomorphisms of $P$ such that $x(t) = \p_{H,t}(x)$ is the solution of 
(\ref{hami}).
The set ${\cal D}_{\o} = 
\{ \p_{H,1} | H \in C^{\infty}(P \x \BR, 
\BR)\}$ of all diffeomorphisms arising in this way is a
subgroup of the group of symplectic diffeomorphisms.
An element in the set ${\cal D}_{\o}$ of 
exact diffeomorphisms is called a (time--dependent) {\em 
exact isotopy.}

For an exact isotopy $\p = \{ \p_t\}_{0 \leq t \leq 1}$
on $(P, \o)$, we define the space
\[\O_{\p} = \{ z\co I \to P \ \ \,  | \ \ z(0) \in L, z(1) \in \p_1(L), 
[\p_t^{-1}z(t) ] = 0 \in \pi_1(P, L) \} . \]
Let $L$ be a monotone Lagrangian sub--manifold,
and let $\p = \{ \p_t\}_{0 \leq t \leq 1}$
be an exact isotopy on $(P, \o)$. If $u$ and $v$ are two maps from 
$[0,1] \x [0,1]$ to $\O_{\p}$ such that
\[u(\t , 0) , v(\t , 0) \in L, \ \ \ u(\t, 1 ) , v(\t, 1) \in \p_1(L) \]
\[u(0,t) = v(0,t) \equiv x, \ \ \ u(1,t) = v(1,t) \equiv y, \ \ \ 
x, y \in L \cap \p_1(L),\]
then we have
\[I_{\o}(u) = I_{\o}(v) \ \ \mbox{if and only if} \ \ \ 
\mu_u(x,y) = \mu_v(x,y), \]
where $\mu_u(x, y) = I_{\mu, L}(u)$ is the Maslov--Viterbo index.
In particular, if $u$ and $v$ are $J$--holomorphic curves with respect to the 
almost complex structure $J$ (may vary with time $t$) compatible with $\o$,
then 
\[ \int \|\n u \|^2_J = \int \| \n v\|^2_J \ \ \ \mbox{if and only if}
\ \ \ I_{\mu, L}(u) = I_{\mu, L}(v) . \]

Note that $\mu_u(x,y)$ is well--defined
mod $\Si (L)$.
The tangent space $T_z\O_{\p}$ of $\O_{\p}$
consists of vector fields $\xi $ of $P$ along 
$z$ which are tangent to $L$ at $0$ and to $\p_1(L)$ at $1$. Then $\o$ induces
a ``1--form'' on $\O_{\p}$:
\begin{equation} \label{1form}
 Da(z) \xi = \int_0^1 \o (\frac{dz}{dt}, \xi(t)) dt .
\end{equation} 
This form is closed in the sense that it can be integrated locally to a real
function $a$ on $\O_{\p}$.
The term $Da(z) \xi$ vanishes
for all $\xi $ if and only if $z$ is a constant loop, i.e., $z(0)$ is 
a fixed point of $\p_1$. The critical point set $Z_{\p }$ 
of the 1--form $Da$ is 
the intersection point set 
$L \cap \p_1(L)$. A critical point is non--degenerate if and only if 
the corresponding intersection is transversal.

For a monotone Lagrangian sub--manifold $L$, an exact isotopy $\p$ and 
$k > 2/p$, consider the space of $L^p_k$--paths
\[{\P }^p_{k,loc}(L, \p;P) = \{ u \in L^p_{k, loc}(\T, P) \ \ 
| \ \ u({\BR } \x \{0\}) \subset L, \ \  
u({\BR } \x \{1\} ) \subset \p(L) \}, \]
where ${\T }= {\BR }\x [0,1] ={\BR } \x iI \subset \BC$. 
Let $S_{\o}$ be the bundle of all 
$J \in \text{End}(TP)$ whose fiber is given by 
\[S_x = \{ J \in \text{End}(T_xP) \ \ | \ \ J^2 = - Id \ \ \mbox{and} \ \
\o(\cdot , J \cdot ) \ \ \mbox{is a Riemannian metric}\}. \]
Let ${\J }=
C^{\infty }([0,1] \x S_{\o})$ be 
the set of time--dependent almost complex structures. Define
\begin{equation} \label{holo}
\pj u(\t, t) = \frac{\bd u(\t, t)}{\bd \t} + J_t \frac{\bd u(\t,t)}{\bd t}
\end{equation}
on ${\P }^p_{k, loc}(L, \p;P)$. 
Then the equation $\pj u = 0$ is translational
invariant in the variable $\t$. 
Let ${\CM }$ be the moduli space ${\CM }_J(L, \p) = \{u\co 
{\BR} \to \in \O_{\p}| \int_{{\BR} \x I} | \frac{\bd u}{\bd \t}|^2$ $< \infty, 
\pj u = 0\}$ of finite actions, and
${\CM }_J(x, y) = 
\{ u \in {\CM }| \ \ \lim_{\t \to + \infty }u = x, \lim_{\t \to
- \infty }u = y; \ \ x, y \in Z_{\p} \}$.
So the moduli space ${\CM }$ is the union 
of $J$-holomorphic curves
$\bigcup_{x, y \in L \cap \p_1(L)} {\CM }_J(x, y)$. If $L$ intersects $\p_1(L)$
transversely, then there exists a smooth Banach
manifold ${\P }(x, y) = {\P }^p_k(x,y) \subset {\P }^p_{k,loc}$ 
for each $x, y \in Z_{\p}$ such that 
(\ref{holo}) defines a smooth section $\pj$ of the smooth Banach space
bundle ${\cal L}$ over ${\P }(x,y)$ 
with fibers ${\cal L}_u = L^p_{k-1}(u^*TP)$.
So ${\CM }_J(x,y)$ is the zero set of $\pj $. The tangent space $T_u{\P}$
consists of all elements $\xi \in L^p_k(u^*(TP))$ so that $\xi (\t, 0) \in TL$ 
and $\xi(\t, 1) \in T(\p_1 (L))$ for all $\t \in \BR$. 
The linearized operator of
$\pj u$, denoted by
\begin{equation*}
E_u = D \pj (u)\co T_u {\P } \to {\cal L}_u,
\end{equation*}
is a Fredholm operator for $u \in {\CM }_J(x,y)$. There is a dense set 
${\J }_{reg}(L,
\p_1(L)) \subset {\J }$ so that if $J \in {\J }_{reg}(L,\p_1(L)) $, 
then $E_u$ is 
surjective for all $u \in {\CM }_J(x, y)$. Moreover the Fredholm index of the 
linearized operator $E_u$ is the same as the Maslov index $\mu_u(x,y)$.
In particular, the space ${\CM }_J(x,y)$ is 
a smooth manifold with dimension 
$\mu_u(x,y)$ for $J \in {\J }_{reg}(L, \p_1(L))$
(see Proposition 2.1 in \cite{fl2}).

\begin{thm}\label{focohomology}{\rm\cite{fl2,oh}}\qua 
Let $L$ be a monotone Lagrangian sub--manifold in $P$, $\Si (L) \geq 3$
 and $\p = \{\p_t\}_{0 \leq
t \leq 1}$ be an exact isotopy such that $L$ intersects 
$\p_1(L)$ transversely. 
Then there is a dense subset
${\J }_*(L, \p) \subset {\J }_{reg} (L, \p)$ of $\J$ such that (1) the zero
dimensional component of ${\hM }_J(x, y) = {\CM}_J(x, y)/{\BR}$
is compact and (2) the one dimensional
component of ${\hM }_J(x^{'}, y^{'})={\CM}_J(x^{'}, y^{'})/{\BR}$ 
is compact up to the splitting of two
isolated trajectories for $J \in {\J }_*(L, \p)$.
Let $C_*(L, \p; P, J)$ be the 
free module over ${\BZ}_2$ generated by $Z_{\p}$.
Moreover, there exists a homomorphism
\begin{equation} \label{25}
 \d \co C_*(L, \p; P, J) \to C_*(L, \p; P, J)
\end{equation}
with $\d \circ \d = 0$ for $J \in {\J }_{reg}(L, \p_1(L))$. The 
${\BZ }_{\Si(L)}$--graded symplectic Floer cohomology
$HF^*(L, \p; P, J)$ is defined to be
the cohomology of the complex $(C_*(L, \p; P, J), \d )$, and
$HF^*(L, \p; P, J)$ is invariant under the continuation of $(J, \p)$,
denoted by $HF^{*}(L, \p; P)$ with $* \in {\BZ}_{\Si (L)}$.
\end{thm}

In order to extend the ${\BZ }_{\Si(L)}$--graded symplectic
Floer--Oh cohomology to a $\BZ$--graded symplectic Floer cohomology,
we make 
use of the infinite cyclic cover $\tilde{\O}_{\p}^*$ of $\O_{\p}$.
By (\ref{1form}), the functional $a$ on
${\O }_{\p}$ is only defined as
$a \co {\O }_{\p} \to {\BR }/ \s (L) \BZ$, for different topology classes in
$\pi_2(P, L)$.
The symplectic
action on $\tilde{\O}_{\p}^*$ and the Maslov index function on 
$\tilde{Z}_{\p}^*$ are well--defined: 
$a\co \tilde{\O}_{\p}^* \to \BR$ and $\mu \co \tilde{Z}_{\p}^* \to \BZ$. 
The ${\BZ}$--graded symplectic Floer cohomology is constructed from the lifted
symplectic action and the lifted Maslov index.
The functional $a$ on ${\O }_{\p}$ and its lift on 
$\tilde{\O}_{\p}^*$ are clearly distinguished from the context. 

\begin{lm} \label{univ}
There exists a universal covering
 space $\tilde{\O}_{\p}$ of ${\O }_{\p}$ with transformation
group $\pi_2 (P, L)$.
\end{lm} 
\proof
By theorem 3.1 
in \cite{mi}, there is a universal covering space 
$\tilde{\O}_{\p}$ of ${\O }_{\p}$ since the space
${\O }_{\p}$  has the homotopy type of a CW complex.
For $[u]\in \pi_1({\O }_{\p})$, we have a representative $u\co I \to {\O }_{\p}$
such that there is a homotopy $F(\t, t)$
of ${\p }^{-1}_t u(\t,t)$ to a constant
path in $(P, L)$ by the definition of ${\O }_{\p}$.
Thus we can reformulate the map
$u$ to yield a map $\overline{u}(\t, t) = u(2 \t, t)$ for $0 \leq \t \leq 1/2$;
$\overline{u}(\t, t) = F(2 \t -1, t)$ for $1/2 \leq \t \leq 1$. Such a map
$\overline{u}\co (D^2, \bd D^2) \to (P, L)$ defines an element in $\pi_2(P,L)$.
It is easy to check that $u \to \overline{u}$ is a bijective homomorphism
between $\pi_1(\O_{\p})$ and $\pi_2(P,L)$ (see also Proposition 2.3
in \cite{fl2}). So the result follows. \endproof

Now the closed 1--form $Da(z)$ has a function 
$a\co \tilde{\O}_{\p} \to \BR$ which is well--defined up to a constant. 
Pick a point $z_0 \in L \cap 
\p_1(L)$ such that $a(z_0) = 0$ by adding a constant.
 For $g 
\in \pi_1(\O_{\p}) = \pi_2(P, L)$, we have 
\begin{equation} \label{deg}
a(g(x)) = a(x)
+ \deg (g) \s (L) , 
\end{equation}
where $\deg (g)$ is defined by $I_{\o}(g) = \deg (g) \s (L)$.
Let $\text{Im}\, (a)(Z_{\p})$ be the image of
$a$ of $Z_{\p}$; modulo $\s(L) \BZ$, the set 
$\text{Im}\, (a)(Z_{\p})$ is finite. Thus
the set ${\BR }_{L, \p} = {\BR } \setminus \text{Im}\, (a)(Z_{\p})$ 
consists of the regular
values of the symplectic action $a$ on $\tilde{\O }_{\p}$. 
From the map $a\co {\O}_{\p} \to {\BR}/{\s (L) {\BZ}}$, we pullback
the universal covering space ${\BR} \to {\BR}/{\s (L) {\BZ}}$ over
${\O}_{\p}$. Let $\tilde{\O}_{\p}^*$ be the pullback
$a^*({\BR}) \to {\O}_{\p}$. The space $\tilde{\O}_{\p}^*$
is an infinite cyclic sub--covering space of the
covering space $\tilde{\O}_{\p}$ of ${\O}_{\p}$.
 
Given $x \in Z_{\p} \subset \O_{\p}$, let $x^{(r)} \in \tilde{Z}^*_{\p} \subset
\tilde{\O}_{\p}^*$ be the unique lift of $x$ such that $a(x^{(r)}) \in
(r, r + \s(L))$. Let $\mu^{(r)}(x) = \mu_{\tilde{u}}(x^{(r)}, z_0) \in \BZ$
for $\tilde{u}= g \circ u$ with $g \in \pi_2(P, L)$ and $x^{(r)}= g(x)$. 
We define
the ${\BZ}$--graded symplectic Floer cochain group by
\begin{equation} \label{cohain}
 C^{(r)}_n(L, \p; P, J) = {\BZ }_2 \{ x \in Z_{\p} \ \ \ | \ \ \ 
\mu^{(r)}(x) = n \in {\BZ  }\} .   
\end{equation}
The group $C^{(r)}_n(L, \p; P, J)$ is a free module over ${\BZ}_2$
generated by $x \in Z_{\p}$ with 
$\mu(x^{(r)}, z_0) = n$. The grading is independent of the choice of $z_0$.
If $\overline{z}_0$ is another choice of the based point and
$g(z_0) = \overline{z}_0$ for some covering transformation $g$, then the
corresponding choice of a lift $\overline{x}^{(r)}$ of $x$ is just $g(x^{(r)})$
by (\ref{deg}).
Note that the Maslov index $\mu^{(r)}(x)$ is independent of the
choice of the based point $z_0$ used in the definition of $a$.
The following 
lemma shows that the lift of the functional $a$ is compatible with 
the universal 
lift of the circle ${\BR }/ \s (L) \BZ$.
\begin{lm} \label{compatible}
The lift of the symplectic action over $\tilde{\O}_{\p}^*$ is 
compatible with the one of the Maslov index: for $g \in \pi_2(P, L)$
with $\deg \, (g) = n$,
\[ a(g(z_0)) = n \s (L) \ \ \ \mbox{if and only if} \ \ \ \ 
\mu^{(r)}(g(z_0), z_0) = n \Si (L). \]
\end{lm} 
\proof  Let  $J$ be a compatible almost complex structure and
$\o (\cdot , J \cdot )$ be the corresponding Riemannian metric on $P$. Denote
$\n $ be the Levi--Civita connection of the metric $\o (\cdot , J \cdot )$.
Then $T_xL$ is an orthogonal complement of $JT_xL$. One can represent
$J_x$ to be the standard $J$ for suitable orthonormal basis in $T_xL$. Let $h$
be the parallel transport along the path $u(\t )$ ($u(\t, t)$ for each fixed
$t \in I$) in $\O_{\p}$. Then we get an isometry
\[h_{\t, t}\co T_xP \to T_{u(\t, t)}P. \]
Define $J_{\t, t} = h_{\t, t}^{-1} \circ J_{u(\t, t)} \circ h_{\t, t}$. Then
we have a smooth map $f\co I \x I \to SO(2m)$ such that $f_{\t, t}^{-1} \circ
J_{\t, t} \circ f_{\t, t} = J_x$. Set
\[ \tilde{L}(\t) = h_{\t, 0}^{-1}(T_{u(\t,0)}L), \ \ \
\widetilde{\p_1(L)}(\t) = h_{\t, 1}^{-1}(T_{u(\t,1)}\p_1(L)). \]
Thus $f_{\t, 0}(\tilde{L}(\t)) = L(\t)$ and $f_{\t, 1}(\tilde{\p_1(L)}(\t)) =
\p_1(L)(\t)$. The trivialization of $u^*TP$ by using the parallel transportation
$\{h_{\t, t}\}$ is given by \[u^*(T_xP) = I \x I \x T_xP = I \x I \x {\BC}^m.\]
Then there are two paths of Lagrangian subspaces $\tilde{L}(\t)$ and 
$\widetilde{\p_1(L)}(\t)$ in
$T_xP = {\BC}^m$. Note that these two Lagrangian paths intersect transversely 
at end points $\t = 0$ and $\t = 1$. There is a map $f_{\p_1}$
from the space $\O_{\p}$ to the space $\Lam (m)$ of pairs
of Lagrangian subspaces in ${\BC}^m$ 
defined by 
$$f_{\p_1}(\{u(\t )\}) = \{L(\t), \p_1 (L)(\t )\}, 0 \leq
\t \leq 1. $$
The Lagrangian Grassmannian $\Lam (m)$ 
has a universal covering $\tilde{\Lam }(m)$  
\cite{clm}.  
For the map $f_{\p_1}\co \O_{\p} \to \Lam (m)$, there is a map from the
CW complex $\O_{\p}$ to $\Lam (m)$ from the obstruction theory. 
Hence there exists
a corresponding map $F$ between the covering space $\tilde{\O}_{\p}^*$
and the universal covering space $\tilde{\Lam }(m)$. From the choice of $z_0$,
$a(g(z_0)) = n \s (L)$. Note that $u_g(0) = z_0, u_g(1) = g (z_0)$,
and $\{u_g(\t )\}_{0 \leq \t \leq 1}$ corresponds to an element
$g \in \pi_2 (P, L)$.

By the definitions of $\s (L)$ and $\Si (L)$, we have $\deg \co
\pi_1 (\O_{\p}) \to \s (L) \BZ$ and 
$Mas \co \pi_1 (\Lam (m)) \cong \Si (L) \BZ$.
So there is a $g_1 \in \pi_1 (\Lam (m))$ induced by $g$ such that 
$g_1 \circ F = f_{\p_1} \circ g$. 
Note that $\mu^{(r)}(z_0, z_0) = 0$. The following
diagram is commutative:
\begin{equation*} 
\begin{CD}
\pi_1 (\O_{\p })@>{\pi_1(f_{\p })}>>\pi_1 ({\Lam } (m))\\
@VV{\deg(g)}V    @VV{\deg(g_1)}V\\
\s (L) {\BZ } @>{F_*}>>\Si (L) \BZ . 
\end{CD}  
\end{equation*}
So $I_{\o }(u_g) = n \s (L)$ and $I_{\mu , L}(u_g) = \deg \, (g_1) \Si (L)$ by
the definitions of $\Si (L)$ and the Maslov index. 
Thus the result follows from $\s (L) = \lam \Si (L)$ and the
monotonicity of $L$.
\endproof
 
The index $\mu_u (x)$ depends on the trivialization
over $I \x I$, only the relative index does not depend on the trivialization.
So the choice of a single $z_0$ fixes the shifting in the $\BZ$--graded 
symplectic Floer
cochain complex. For $g\in \pi_2(P, L)$ with $x^{(r)}=g(x)$, we have
$\mu^{(r)}(x) = \mu (x) + \deg g \cdot \Si (L)$.

\begin{pro}[Lemma~\ref{compatible} and Proposition 2.4 \cite{fl2}]\label{bub}
If $u \in {\P }(x, y)$ for $x, y \in Z_{\p}$ and $\tilde{u}$ is any
lift of $u$, then $\mu_{\tilde{u}}(x^{(r)}, y^{(r)}) 
= \mu^{(r)}(y) - \mu^{(r)}(x) = \mu (y^{(r)}, z_0) - \mu (x^{(r)}, z_0)$.
\end{pro}
 
\begin{df} \label{integ} 
The ${\BZ}$--graded symplectic Floer coboundary map 
is defined by
$$\bd^{(r)}\co C^{(r)}_{n-1}(L, \p; P, J) \to C^{(r)}_{n}(L, \p; P, J)$$
\[ \bd^{(r)}x = \sum_{y \in C^{(r)}_{n}(L, \p; P, J)} \# {\hM}_J(x, y) \cdot
 y, \]
where ${\CM }_J(x,y)$ is the union of the components of 1--dimensional moduli 
space of $J$--holomorphic curves,
and ${\hM }_J(x,y) = {\CM }_J(x, y)/{\BR }$ is the 
zero--dimensional moduli space modulo $\t$--translational invariant. 
The number $\# {\hM }_J(x, y)$
counts the points modulo 2.
\end{df} 

\noindent{\bf Remark}\qua The condition $\Si (L) \geq 3$, rather
than $\Si (L) \geq 2$, enters only in proving that 
$\la \d \circ \d x, x\ra = 0$.
For $\Si (L) = 2$, Oh 
evaluated a number $\pmod 2$ of $J$--holomorphic
disks with Maslov index $2$ that pass through
$x \in L \subset P$, and
verified that the number is always even. Hence
$\la \d \circ \d x, x\ra = 0 \pmod 2$. In our case, this reflects to understand
the two lifts $x^{(r)}$ and $g(x^{(r)})$ of $x$ with $\deg \, (g) = \pm 1$.
Note that $x^{(r)} \in (r , r + \s (L))$ and $g(x^{(r)}) \in
(r + \deg (g) \s (L), r+ (\deg (g) +1) \s (L))$. So
the ${\BZ}$--graded symplectic
coboundary is not well--defined in this case. We leave it to future study.

The coboundary map $\bd^{(r)}$ only counts part of Floer's 
coboundary map in (\ref{25}). Next task is to verify 
$\bd^{(r)} \circ \bd^{(r)} = 0$ in the following.

\begin{lm} Under the same hypothesis in Theorem~\ref{focohomology}, we have
$\bd^{(r)} \circ \bd^{(r)}$ $= 0$.
\end{lm}
\proof 
If $x \in C^{(r)}_{n-1}(L, \p; P, J)$ ($\mu (x^{(r)}, z_0) = n-1$), 
by the definition of $\bd^{(r)}$,
then the coefficient of $z \in C^{(r)}_{n+1}(L, \p; P, J)$  in 
$\bd^{(r)} \circ \bd^{(r)}(x)$ is given by
\begin{equation}
 \sum_{y \in C^{(r)}_{n}(L, \p; P, J)} 
\# {\hM }_J(x, y) \cdot \# {\hM }_J(y,z). 
\end{equation}
By Proposition~\ref{bub}, the boundary of the 1--dimensional
manifold ${\hM}_J(x, z) = {\CM }_J(x, z)/{\BR}$ corresponds to two isolated
trajectories
${\CM }_J(x, y) \x {\CM }_J(y, z)$. 
Each term $\# {\hM}_J(x, y) \cdot \# {\hM}_J(y,z)$ is
the number of the 2--cusp trajectories of ${\hM}_J(x, z)$ with 
$y \in C^{(r)}_{n}(L, \p; P, J)$. 
For any such $y$, there are $J$--holomorphic curves $u \in {\CM }_J(x, y)$ and
$v \in {\CM }_J(y,z)$.
The other end of the corresponding 
component of ${\hM}_J(x,z)$ corresponds to the space
${\CM }_J(x, y^{'}) \x {\CM }_J(y{'},z)$ with 
$u^{'} \in {\CM }_J(x, y^{'})$ and 
$v^{'} \in {\CM }_J(y{'},z)$. Then ${\hM}_J(x, z)$ has an 
1--parameter family of paths 
from $x$ to $z$ with ends $u \# v$ and $u^{'} \# v^{'}$ for appropriate
grafting (see \cite{fl2} section 4). 
If we lift $u$ to $\tilde{u} \in \tilde{\CM}_J(x^{(r)}, 
\tilde{y})$ the moduli space of J--holomorphic curves in $\tilde{\O}_{\p}$
with asymptotics $x^{(r)}$ and $\tilde{y}$, then
\begin{equation} \label{rest}
1 = I_{\mu, \tilde{L}}(\tilde{u})=
\mu(\tilde{y}, z_0) - \mu^{(r)}(x) = 
\mu(\tilde{y}, z_0) - (n-1).
\end{equation}
So $\mu(\tilde{y}, z_0) = n$, and $\tilde{y} = y^{(r)}$ is the preferred
lift. So $\mu^{(r)}(y) = \mu (\tilde{y}, z_0) = n$.
Thus $\tilde{u} \in \tilde{\CM}_J(x^{(r)}, 
y^{(r)})$. Similarly $\tilde{v} \in \tilde{\CM}_J(y^{(r)}, z^{(r)})$.
Since $u^{'} \# v^{'}$ is homotopic to $u \# v$ rel $(x^{(r)}, z^{(r)})$, the 
lift $\tilde{u}^{'} \# \tilde{v}^{'}$ is also a path with ends 
$(x^{(r)}, z^{(r)})$. Using the fact of the symplectic action $a$ is 
non-increasing along any gradient trajectory $\tilde{u}^{'}$, we have
\begin{equation}
r < a(z^{(r)}) \leq a(\tilde{y}^{'}) \leq a(x^{(r)}) < r + \s (L).
\end{equation} \label{ineq}
By the uniqueness, we have $\tilde{y}^{'} = (y^{'})^{(r)}$.
By (\ref{rest})
for $u^{'}$, we have $\mu^{(r)}((y^{'})^{(r)}) = \mu^{(r)}(x^{(r)}) +1 = n$.
So $y^{'} \in C^{(r)}_{n}(L, \p; P, J)$. Thus the number of
two--cusp trajectories connecting $x^{(r)}$ and $z^{(r)}$ 
with index 2 is always even. Hence we obtain $\bd^{(r)} \circ \bd^{(r)} = 0$. 
\endproof

The complex $(C^{(r)}_{n}(L, \p; P, J), \bd^{(r)}_n)_{n\in \BZ}$ 
is indeed a 
$\BZ$--graded symplectic Floer cochain complex. 
We call its cohomology to be an {\em 
${\BZ}$--graded symplectic Floer cohomology},
denoted by
\begin{equation} \label{zgs}
I^{(r)}_*(L, \p; P, J) = H^*(C^{(r)}_*(L, \p; P, J), \bd^{(r)}), \ \ \
* \in \BZ. \end{equation}
By the construction of $I^{(r)}_*(L, \p; P, J)$,
if $[r,s] \subset {\BR }_{L, \p}$, then
$I^{(r)}_*(L, \p; P, J) = I^{(s)}_*(L, \p; P, J)$.
The relation between $I^{(r)}_*(L, \p; P, J)$ and $HF^*(L, \p; P)$ will be 
discussed in \S 4.

\section{Invariance property of the ${\BZ}$--graded symplectic Floer cohomology}

In this section we show that the $\BZ$--graded symplectic Floer 
cohomology in (\ref{zgs})
is invariant under the changes of
$J$ and under the exact deformations $\p_1$ of the Lagrangian sub--manifold $L$.

Let $\{(J^{\lam }, \p^{\lam })\}_{\lam \in \BR}$ be an 1--parameter family 
which interpolates from $(J^0, \p^0)$ to $(J^1, \p^1)$.
The family $(J^{\lam }, \p^{\lam })$ is constant in
$\lam $ outside $[0,1]$. 
We also assume
that $\p_1^{\lam }$ is exact under the change of $\lam $. Let $J^{\lam }_t 
= J(\lam , t)$ be a 2-parameter family 
of almost complex structures compatible to $\o $, and 
$\p^{\lam }_t = \p (\lam, t)$ with $\p (\lam , 0) = \text{Id}$ 
is the 2--parameter 
family of exact isotopies contractible to the identity. 
Such $\p^{\lam }_t$ connecting 
$\p_t^0, \p_t^1$ does exist.
Floer \cite{fl3} discussed the invariance of the symplectic Floer cohomology
under the change of $(J, \p)$ for
$(J^0, \p^0)$ $C^{\infty }$--close to $(J^1, \p^1)$.
Let $H^{\lam }$ be the Hamiltonian function generated by
$\p^{\lam } = \{\p_t^{\lam }\}_{0 \leq t \leq 1}$. 
Then the deformed gradient flow of $a_H$ is
\begin{equation} \label{jl}
{\pj }_{\lam } u_{\lam } (\t, t ) + \n_J H_t^{\lam }(u_{\lam } (\t, t )) =
\frac{\bd u_{\lam }}{\bd \t} + J^{\lam}_t
\frac{\bd u_{\lam }}{\bd t} + \n_J H_t(u_{\lam } (\t, t )) = 0 ,
\end{equation}
with the moving Lagrangian boundary conditions
\begin{equation} \label{ll}
u_{\lam }(\t, 0) \in L,  u_{\lam }(\t, 1) \in \p_1^{\lam}(L).
\end{equation}
We define
$C^{(r)}_{J, L} = \min \{ a(x^{(r)}) -r, \s (L) +r - a(x^{(r)}) | x \in
Z_{\p } \}$.
For each $x \in Z_{\p}$, there is an open neighborhood $U_x$ in ${\O}_{\p}$
such that (1) $U_x$ is evenly covered in
$\tilde{{\O}}_{\p}^*$, (2) for each $z \in U_x$,
$|a(z) - a(x)| < C^{(r)}_{J, L}/8$.
There are finite sub--cover $\{U_{x_1},
\cdots, U_{x_k}\}$ of $Z_{\p}$, and by Gromov's compactness theorem
\cite{fl3} we have $\ve_1 > 0$ such that if $\|Da(z)\|_{L^3_1} < \ve_1$ then
$z \in \bigcup_{i=1}^kU_{x_i}$. Let $\ve = \min \{ \ve_1,
C^{(r)}_{L_1, L_2}/8\}$. We set a deformation $\{J, \p\}$
satisfying the usual perturbation requirements in \cite{fl3}, and also
satisfying
\begin{equation} \label{admi}
{\rm(i)} \ \ |H_t^{\lam}(z)| < \ve /2, \ \ \ \
{\rm(ii)} \ \ \|\n_J H_t^{\lam}(z)\|_{L^3_1} < \ve /2,
\end{equation}
for all $z \in {\O}_{\p}$. These deformation conditions can be achieved by
the density statement in \cite{fl3}.

Let ${\P }_{1, \ve /2}$
be the set of $\{J, \p \}$ which satisfies these extra conditions (\ref{admi}).

This directly generalizes the $J$-holomorphic curve equation in the cases
of $(J^0, \p^0)$ and $(J^1, \p^1)$. The moduli space ${\CM }_{J^{\lam}}(x, y)$
of (\ref{jl}) and (\ref{ll}) has the same analytic properties as the 
moduli space ${\CM }_J(x, y)$ except for the translational invariance (see
Proposition 3.2 in \cite{fl2}). Hofer analyzed the compactness 
property for a similar moving Lagrangian coboundary condition, Oh \cite{oh}
determined that the bubbling--off spheres or disks 
can not occur in the components
of ${\CM }_{J^{\lam}}(x, y)$ for the monotone Lagrangian sub--manifold $L$ with 
${\Si }(L) \geq 3$. The index of $u_{\lam }$ can be proved to be the same 
as a topological index for the moduli space 
of perturbed $J$--holomorphic curves.
The proof of the invariance under the changes of $(J, \p)$ is 
the same as in \cite{fl2, fl3, oh}. It is sufficient for us
to verify that
the cochain map is well-defined for the $\BZ$--graded symplectic Floer
cochain complexes.
\begin{lm} \label{lift}
 If $u_{\lam} \in {\CM }_{J^{\lam}} (x_0, x_1)$, 
$(J^{\lam }_t, \p^{\lam }_t) \in {\P }_{1, \ve /2}$, and 
$\tilde{u}_{\lam } \in {\P }(\tilde{x}_0, \tilde{x}_1)$ 
is any lift of $u_{\lam }$, then
\[a_{(J^1, \p^1)}(\tilde{x}_1) < a_{(J^0, \p^0)} (\tilde{x}_0) + \ve . \]
\end{lm} 
\proof 
 Note that the path $\{u_{\lam }(\t ) | \t \in (- \infty , 0) \}$ is
a gradient trajectory for $(J^0, \p^0)$ and $\{u_{\lam }(\t ) | \t \in 
(1, \infty )\}$ is a gradient trajectory for $(J^1, \p^1)$. So 
\begin{equation} \label{ine}
 a_{(J^0, \p^0)} (\tilde{u}(0)) \leq a_{(J^0, \p^0)} (\tilde{x}_0), 
\ \ \ a_{(J^1, \p^1)}(\tilde{x}_1) \leq a_{(J^1, \p^1)}(\tilde{u}(1)) .
\end{equation} 
Since $u_{\lam } \in {\CM }_{J^{\lam}}(x_0, x_1)$, by the property of 
${\P }_{1, \ve /2}$, we have
\begin{eqnarray*}
\|{\pj }_{\lam } u_{\lam } (\t, t )\|_{L_1^3} & = &
\|\frac{\bd u_{\lam }}{\bd \t} + J^{\lam}_t
\frac{\bd u_{\lam }}{\bd t}\|_{L_1^3} \\
& = & \|\n_J H_t^{\lam}(u_{\lam } (\t, t ))\|_{L_1^3}\\
& < & \ve /2.
\end{eqnarray*}
By the Sobolev embedding $L_1^3 \hookrightarrow L^2$,
\begin{eqnarray*}
I_{\o }({u}_{\lam})|_{\T \x [0,1]} & = &
\|\bd_{J^{\lam }} {u}_{\lam }\|^2_{L^2(\T \x [0,1])} -
\|{\pj }_{\lam }{u}_{\lam}\|^2_{L^2(\T \x [0,1])} \\
 & \geq & - \|{\pj }_{\lam }{u}_{\lam}\|^2_{L^2(\T \x [0,1])} \\
 & \geq & - \ve/2
\end{eqnarray*}
Thus we obtain the following.
\begin{eqnarray} \label{nedd}
a_{(J^0, \p^0)}(\tilde{u}_{\lam}(0)) & = & a_{(J^1, \p^1)}(\tilde{x}_1) +
I_{\o }(u_{\lam})|_{\T \x [0,1]} \\
& > & a_{(J^1, \p^1)}(\tilde{x}_1)  - \ve/2. \nonumber
\end{eqnarray}
\begin{eqnarray*}
a_{(J^0, \p^0)}(\tilde{x}_0) & \geq & a_{(J^0, \p^0)}(u_{\lam}(0)) \\
& = & a_{(J^0, \p^0)}(u_{\lam}(0)) + H_0(u_{\lam}(0))\\
& > & (a_{(J^1, \p^1)}(\tilde{x}_1)  - \ve/2) - \ve/2,
\end{eqnarray*}
by (\ref{admi}) and (\ref{nedd}).
Then the result follows. \endproof

\begin{df} \label{cochain}
For $n \in \BZ$, define a homomorphism
$\p_{01}^n \co C^{(r)}_n (L, \p^0; P, J^0) \to C^{(r)}_n (L, \p^1; P, J^1)$ by
\[ \p_{01}^n (x_0) = \sum_{x_1 \in C^{(r)}_n (L, \p^1; P, J^1)} \#
{\CM }^0_{J^{\lam}}(x_0, x_1) \cdot x_1, \]
where ${\CM }^0_{J^{\lam}}(x_0, x_1)$ is a 0--dimensional moduli space
of $J$--holomorphic curves satisfying (\ref{jl}) and (\ref{ll}).
\end{df}

We show that the homomorphism $\p_{01}^n$ is a cochain map with respect
to the integral lifts.

\begin{lm} \label{01map}
The homomorphism $\{\p_{01}^*\}_{* \in Z}$ is a cochain map: 
\[ \bd^{(r)}_{n,1} \circ \p_{01}^n = \p_{01}^n \circ \bd^{(r)}_{n,0}, 
\ \ \ \  n \in {\BZ} . \]
\end{lm}
\proof  For $x_0 \in C^{(r)}_n (L, \p^0; P, J^0)$ and 
$y_1 \in C^{(r)}_{n+1} (L, \p^1; P, J^1)$, the coefficient of $y_1$ in
$(\bd^{(r)}_{n,1} \circ \p_{01}^n  - \p_{01}^n \circ \bd^{(r)}_{n,0})(x_0)$
is the modulo 2 number of the set:
\begin{equation} \label{42}
 \bigcup_{y_0 \in C^{(r)}_{n+1} (L, \p^0; P, J^0)} \hat{\CM }_{J^0}
(x_0, y_0) \x {\CM }_{J^{\lam}}^0 (y_0, y_1) 
\end{equation}
\[\coprod \bigcup_{x_1 \in
C^{(r)}_n (L, \p^1; P, J^1)} {\CM }_{J^{\lam}}^0 (x_0, x_1) \x
\hat{\CM }_{J^1}(x_1, y_1).\]
The ends of the 1--dimensional manifold
${\CM }_{J^{\lam}}^1(x_0, y_1)$ are in one--to--one correspondence with the set
\begin{equation} \label{43}
 (\bigcup_{y \in Z_{{\p }^0 }} {\hM }_{J^0}
(x_0, y) \x {\CM }_{J^{\lam}}^0(y, y_1)) \cup (\bigcup_{x \in Z_{{\p }^1}}
{\CM }_{J^{\lam}}^0 (x_0, x) \x 
{\hM }_{J^1}(x, y_1)). 
\end{equation}
For an end $u \# v$ of ${\CM }_{J^{\lam}}^1(x_0, y_1)$ corresponding to an 
element in (\ref{42}), the other end $u^{'} \# v^{'}$ of the same 
component 
corresponds to an element in (\ref{43}) (see \cite{fl2} section 4 for the
gluing construction on $u\# v$). For $u \in {\CM }^0_{J^{\lam}}(x_0, y)$
and $v \in {\hM }_{J^1}(y, y_1)$, 
the space ${\CM }^1_{J^{\lam}}(x_0, y_1)$
gives a 1-parameter family of paths in ${\P }(L, \p^{\lam }; P)$ with 
fixed end points $x_0$ and $y_1$. The 1--parameter family
 gives the homotopy of paths from
$u \#_{\r } v$ to $u^{'} \#_{\r } v^{'}$ rel end points. The lift of 
$u \#_{\r } v$ starts at $x_0^{(r)}$ and ends at $y_1^{(r)}$, so does
the lift of $u^{'} \#_{\r } v^{'}$. 
Suppose $u^{'}$ lifts to an element in ${\CM }^0_{J^{\lam}}
(x_0^{(r)}, \tilde{y})$.
By Lemma~\ref{lift}, we have 
\begin{equation} \label{pf1}
 a_{(J^1, \p^1)} (\tilde{y}) < a_{(J^0, \p^0)} (x_0^{(r)}) + \ve 
< r + \s (L). \end{equation}
Using the fact that trajectory decreases the symplectic action, we obtain 
\begin{equation} \label{pf2}
a_{(J^1, \p^1)} (\tilde{y}) > a_{(J^1, \p^1)} (y_1^{(r)}) > r. 
\end{equation}
So $a_{(J^1, \p^1)} (\tilde{y}) \in (r, r+\s (L))$. 
Inequalities (\ref{pf1}) and (\ref{pf2}) give the
preferred lift $\tilde{y} = y^{(r)}$ of $y$. 
By Proposition~\ref{bub} (1), 
\[ 1 = \mu^{(r)}(y_1) - \mu^{(r)}(y) 
= (n+1) - \mu^{(r)}(y). \]
So $\mu^{(r)}(y) = n$ and $y \in C_n^{(r)}(L, \p^1; P, J^1)$. This
shows that the 
$u^{'} \#_{\r } v^{'}$ in (\ref{43}) actually corresponds to an 
element in (\ref{42}). So the cardinality is always even. \endproof

For $(J^i, \p^i) \in {\P }_{1, \ve /2}$ $(i = 0, 1, 2)$, we define a 
class ${\P }_{2, \ve}$ of perturbations consisting of 
\begin{equation*}
(J^{\lam }, {\p }^{\lam }) =  
\begin{cases}
(J^0, \p^0) & \text{if $\lam \leq - T$}, \\
(J^1, \p^1) & \text{if $- T +1 \leq \lam \leq T -1$}, \\
(J^2, \p^2) & \text{if $\lam \geq T$},
\end{cases} 
\end{equation*}
for a fixed number $T ( > 2)$ with $(J^{\lam }, {\p }^{\lam })_{
0 \leq \lam \leq 1} \in {\P }_{1, \ve /2}$ and $(J^{\lam }, {\p }^{\lam })_{
1\leq \lam \leq 2} \in {\P }_{1, \ve /2}$.
If both perturbations 
$(J^{\lam }, {\p }^{\lam })_{0\leq \lam \leq 1} \in 
{\P }_{1, \ve /2}((J^0, \p^0), (J^1, \p^1))$
and $(J^{\lam }, {\p }^{\lam })_{1\leq \lam \leq 2} \in 
{\P }_{1, \ve /2}((J^1, \p^1), (J^2, \p^2))$,
then we can compose $(J^{\lam }, {\p }^{\lam })_{0\leq \lam \leq 1}$ with 
$(J^{\lam }, {\p }^{\lam })_{1\leq \lam \leq 2}$ to get 
$(J^{\lam }, {\p }^{\lam }) \in
{\P }_{2, \ve }((J^0, \p^0),(J^2, \p^2))$. 
Let $(J^{\lam }, {\p }^{\lam }) =
(J^{\lam }, {\p }^{\lam })_{0\leq \lam \leq 1} \#_T 
(J^{\lam }, {\p }^{\lam })_{1\leq \lam \leq 2}$ be the composition.
Then for a large fixed $T$ and each compact set $K$ in
${\CM }_{J^{\lam }_{0\leq \lam \leq 1}}(x, y) \x
{\CM }_{J^{\lam }_{1\leq \lam \leq 2}}(y, z)$, 
there is a ${\r }_T > 0$ and 
for all $\r > {\r }_T$ a local diffeomorphism
\begin{equation} \label{pf3}
 \#_{{\r }_T} \co K \to 
{\CM }_{(J^{\lam }, {\p }^{\lam })_{0\leq \lam \leq 1} \#_{T,{\r }_T} 
(J^{\lam }, {\p }^{\lam })_{1\leq \lam \leq 2}}
(x, z). 
\end{equation}
See Proposition 2d.1 in \cite{fl3}. 
\begin{lm} \label{compose}
For $(J^{\lam }, {\p }^{\lam }) \in {\P }_{2, \ve}$, and $\r > {\r }_T$, we have
\[ \p_{02}^n = \p_{12}^n \circ \p_{01}^n ,  \ \ \ \ 
\mbox{for $n \in \BZ$} . \]
\end{lm}
\proof  For $x_0 \in C_n^{(r)}(L, \p^0; P, J^0)$, we have
\[ \p_{02}^n (x_0 ) = \sum_{y_0 \in C_n^{(r)}(L, \p^2; P, J^2)} \# 
{\CM }^0_{J^{\lam}}(x_0, y_0) \cdot y_0, \]
where the summation $\sum $ runs over $y_0 \in C_n^{(r)}(L, \p^2; P, J^2)$.
Also we have 
\[ \p_{12}^n \circ \p_{01}^n (x_0 ) = \sum \# 
({\CM }^0_{J^{\lam }_{0\leq \lam \leq 1}}(x_0, y) \x 
{\CM }_{J^{\lam }_{1\leq \lam \leq 2}}(y,y_0) ) \cdot y_0, \]
where the summation $\sum $ runs over $y \in C_n^{(r)}(L, \p^1; P, J^1)$.
The local diffeomorphism $\#_{{\r }_T}$ in (\ref{pf3}) 
determines the following: 
\[ \# {\CM }^0_{J^{\lam}}(x_0, y_0) = 
\# ({\CM }^0_{J^{\lam }_{0\leq \lam \leq 1}}(x_0, y) \x 
{\CM }^0_{J^{\lam }_{1\leq \lam \leq 2}}(y,y_0)) .\]
All we need to check is that $y \in C_n^{(r)}(L, \p^1; P, J^1)$.
This can be verified by 
the same argument in the Lemma~\ref{01map}. \endproof

For two classes $(J^{\lam }, {\p }^{\lam })$ and $(\overline{J}^{\lam }, 
\overline{{\p }}^{\lam })$ in
${\P }_{2, \ve}((J^0, \p^0),(J^2, \p^2))$, the following lemma shows that the 
induced cochain maps $\p_{02}^n$ and $\overline{\p }_{02}^n$ 
are cochain homotopic
to each other.

\begin{lm} \label{homotopy}
If $(J^{\lam }, {\p }^{\lam })$, $(\overline{J}^{\lam }, 
\overline{{\p }}^{\lam })\in
{\P }_{2, \ve} ((J^0, \p^0),(J^2, \p^2))$ 
can be smoothly deformed from one to another
by a 1-parameter family $(J^{\lam }_s, {\p }^{\lam }_s)$ of $s \in [0, 1]$:
$(J^{\lam }_s, {\p }^{\lam }_s) = (J^{\lam }, {\p }^{\lam })$ for
$s \leq 0$, and $(J^{\lam }_s, {\p }^{\lam }_s) =  (\overline{J}^{\lam }, 
\overline{{\p }}^{\lam })$ for $s \geq 1$. Then $\p_{02}^*$ and
$\overline{\p }_{02}^*$ are cochain homotopic to each other.
\end{lm}
\proof  It suffices to construct a homomorphism
\[ H\co C^{(r)}_*(L, \p^0; P, J^0) \to C^{(r)}_*(L, \p^2; P, J^2) , \]
of degree $- 1$ with the property
\begin{equation} \label{hom}
 \p_{02}^n - \overline{\p }_{02}^n = H {\bd }_{n,0}^{(r)} + 
{\bd }_{n,2}^{(r)} H,  \ \ \ \ \mbox{for $n \in \BZ$}.
\end{equation}
Associated to the family $(J^{\lam }_s, {\p }^{\lam }_s)$, there is 
a moduli space $H{\CM }(x_0, y_0) = \cup_{s \in [0,1]}
{\CM }^0_{(J^{\lam }_s, {\p }^{\lam }_s)}(x_0, y_0)$:
\[ H{\CM }(x_0, y_0) = \{ (u, s) \in 
{\CM }^0_{(J^{\lam }_s, {\p }^{\lam }_s)}(x_0, y_0) \x [0, 1] \} \subset 
{\P }(L, {\p }^{\lam }_s; P)(x_0, y_0) \x [0,1]. \]
The space $H{\CM }(x_0, y_0)$ is the regular zero set of 
${\pj }_{(J^{\lam }_s, {\p }^{\lam }_s)}$, and is smooth manifolds of
dimension $\mu^{(r)}(y_0) - \mu^{(r)}(x_0) + 1$.
For the case of $\mu^{(r)}(x_0) = \mu^{(r)}(y_0) = n$, 
the boundaries of the 1-dimensional manifold $H{\CM }(x_0, y_0)$ of 
${\P }(L, {\p }^{\lam }_s; P)(x_0, y_0) \x [0,1]$ consist of
\begin{itemize}
\item ${\CM }^0_{(J^{\lam }, {\p }^{\lam })}(x_0, y_0) \x \{0 \} 
\cup {\CM }^0_{(\overline{J}^{\lam }, 
\overline{{\p }}^{\lam })}(x_0, y_0) \x \{1 \} $,
\item $\ov{\CM}_{(J^{\lam }_s, {\p }^{\lam }_s)}^0(x_0, y) 
\x {\CM }^0_{(J^2, \p^2)}(y, y_0)$ for $y \in C_{n-1}^{(r)}(L, \p^2; P, J^2)$,
\item ${\CM }^0_{(J^0, \p^0)}(x_0,x) \x 
\ov{\CM}_{(J^{\lam }_s, {\p }^{\lam }_s)}^0(x, y)$ for
$x \in C_{n-1}^{(r)}(L, \p^0; P, J^0)$.
\end{itemize}

Note that $\ov{\CM}_{(J^{\lam }_s, {\p }^{\lam }_s)}^0(x_0, y)$
and $\ov{\CM}_{(J^{\lam }_s, {\p }^{\lam }_s)}^0(x, y)$ are moduli spaces
of solutions of
$(u, s)$ of $J$-holomorphic equations lying in virtual dimension $-1$
($\mu_u = - 1$), they can
only occur  for $0 < s < 1$. Define $H\co C^{(r)}_n(L, \p^0; P, J^0) \to 
C^{(r)}_{n-1}(L, \p^2; P, J^2)$ by
\begin{equation}
H(x_0) = \sum_{ y \in C^{(r)}_{n-1}(L, \p^2; P, J^2)} 
\# \ov{\CM }^0_{(J^{\lam }_s, {\p }^{\lam }_s)}(x_0, y) \cdot y,
\end{equation} 
for $0 < s < 1$.
Similar to Lemma~\ref{01map}, by checking the corresponding preferred lifts and
the integral Maslov indexes, we get the desired cochain homotopy $H$ 
between
${\p }_{02}^n$ and $\overline{\p }_{02}^n$ such that
$H$ satisfies (\ref{hom}). \endproof

By Lemma~\ref{homotopy}, the homomorphism $\p_{02}^*$ induced 
from $(J^{\lam }, 
{\p }^{\lam })$ is the same homomorphism $\overline{\p }_{02}^*$ induced from
$(\overline{J}^{\lam }, \overline{\p }^{\lam })$ on the 
$\BZ$-graded symplectic Floer cohomology. 
So the
$\BZ$-graded symplectic Floer cohomology is invariant
under the continuation of $(J,\p)$. The following is Theorem A (1).

\begin{thm} \label{invariant}
For any continuation $(J^{\lam }, {\p }^{\lam }) \in {\P }_{1, \ve /2}$ which 
is regular at the ends, there exists an isomorphism
\[\p_{01}^n\co I^{(r)}_n(L, \p^0; P, J^0) \to I^{(r)}_n(L, \p^1; P, J^1), 
\ \ \ \ \mbox{for $n \in \BZ$} . \]
\end{thm}
\proof  Let $(J^{- \lam }, {\p }^{ - \lam })$ be the reversed family 
of $(J^{\lam }, {\p }^{\lam })$ by setting $\t = - {\t }^{'}$. So we can form
a family of composition 
$(J^{\lam }, {\p }^{\lam })_{0\leq \lam \leq 1} \#_T 
(J^{- \lam }, {\p }^{ - \lam })_{1\leq \lam \leq 2}$ in 
${\P }_{2, \ve}$ for some fixed $T ( > 2)$. By Lemma~\ref{compose}, 
\[{\p }_{(J^{\lam }, {\p }^{\lam })_{0\leq \lam \leq 1} \#_T 
(J^{- \lam }, {\p }^{ - \lam })_{1\leq \lam \leq 2}} =
{\p }_{10}^* \circ {\p }_{01}^* . \] 
One can deform $(J^{\lam }, {\p }^{\lam })_{0\leq \lam \leq 1} \#_T 
(J^{- \lam }, {\p }^{ - \lam })_{1\leq \lam \leq 2}$
into the trivial continuation $(J^0, \p^0)$ for all $\t \in \BR$.
Then by Lemma~\ref{homotopy}, we have 
\[ {\p }_{10}^* \circ {\p }_{01}^* = {\p }_{00}^* = id:
I^{(r)}_*(L, \p^0; P, J^0) \to I^{(r)}_*(L, \p^0; P, J^0). \]
Similarly, ${\p }_{01}^* \circ {\p }_{10}^* = {\p }_{11}^* = id$ on
$I^{(r)}_*(L, \p^1; P, J^1)$. The result follows. \endproof

The ${\BZ}$--graded symplectic Floer cohomology $I_*^{(r)}$ is functorial 
with respect to compositions of continuations $(J^{\lam }, {\p }^{\lam })$,
and invariant under continuous deformations of $(J^{\lam }, {\p }^{\lam })$
within the set of continuations ${\P }_{1, \ve /2}$.

\section{The spectral sequence for the symplectic Floer cohomology}

In this section, we show that 
$I_*^{(r)}(L, \p;P)$ 
(the $\BZ$--graded symplectic Floer
cohomology) for $r \in {\BR }_{L, \p}$ and $* \in \BZ$ determines
the ${\BZ }_{\Si (L)}$--graded symplectic Floer cohomology 
$HF^*(L, \p; P)$. 
The way to link them together is to filter the 
$\BZ$-graded symplectic Floer cochain complex. Then by a standard method in 
algebraic topology (see \cite{sp}),
the filtration gives arise a spectral sequence
which converges to the ${\BZ }_{\Si (L)}$--graded symplectic Floer cohomology 
$HF^*(L, \p; P)$. The K\"{u}nneth formula (Theorem~\ref{kus}) for the 
spectral sequence is obtained by the analysis of higher differentials and 
Maslov indexes.

\begin{df} \label{filter}
For $r \in {\BR }_{L, \p}, j \in {\BZ }_{\Si (L) }$ and 
$n \equiv j \ \ (\ mod \ \ \Si (L))$, 
we define the free module over ${\BZ}_2$:
\[ F_n^{(r)} C_j(L, \p; P, J) = \sum_{k \geq 0} C^{(r)}_{n + \Si (L) k}
(L, \p; P, J). \]
The free module $F^{(r)}_*C_*(L, \p; P, J)$
gives a natural decreasing filtration on the symplectic Floer
cochain groups $C_*(L, \p; P, J)\, (* \in 
{\BZ }_{\Si (L) })$. 
\end{df}

There is a finite length decreasing filtration of $C_j(L, \p; P, J)$, 
$j \in {\BZ }_{\Si (L) }$:
\begin{equation} \label{dea}
\cdots F_{n+\Si (L)}^{(r)}C_j(L, \p; P, J) \subset F_n^{(r)}
  C_j(L, \p; P, J) \subset 
  \cdots \subset C_j(L, \p; P, J) .
\end{equation}
\begin{equation} \label{coll}
C_j(L, \p; P, J) = \bigcup_{n \equiv j (mod \ \ \Si (L))} 
F_n^{(r)} C_j(L, \p; P, J).
\end{equation}
Note that the symplectic action is non--increasing 
along the gradient trajectories.
The coboundary map
$\d\co F_n^{(r)} C_j(L, \p; P, J) \to F_{n+1}^{(r)} C_{j+1}(L, \p; P, J)$ 
in the Theorem~\ref{focohomology}
preserves the
filtration in Definition~\ref{filter}. Thus the 
${\BZ }_{\Si (L)}$-graded symplectic Floer cochain complex 
$(C_j(L, \p; P, J), \d)_{j \in {\BZ }_{\Si (L)}}$
has a decreasing bounded 
filtration $(F_n^{(r)}C_*(L, \p; P, J), \d)$:
\begin{equation} \label{pre}
\begin{array}{ccccc}
  \downarrow & \downarrow &  \downarrow \\
\cdots  F_{n+\Si (L)}^{(r)}C_j(L, \p; P, J) & 
\subset F_n^{(r)}C_j(L, \p; P, J) & 
\cdots  \subset  C_j(L, \p; P, J) \\
  \downarrow \partial^{(r)} & \downarrow \partial^{(r)} &   \downarrow \d \\
\cdots  F_{n+\Si (L) +1}^{(r)}C_{j+1}(L, \p; P, J) & 
\subset F_{n+1}^{(r)}C_{j+1}(L, \p; P, J)
& \cdots \subset  C_{j+1}(L, \p; P, J) \\
 \downarrow & \downarrow &   \downarrow
\end{array} 
\end{equation}
The cohomology of the vertical cochain subcomplex 
$F_n^{(r)} C_*(L, \p; P, J)$ in the filtration (\ref{pre}) is 
$F_n^{(r)} I_j^{(r)}(L, \p; P, J).$

\begin{lm} There is a filtration for the
$\BZ$--graded symplectic Floer cohomology 
$\{I_*^{(r)}(L, \p; P, J) \}_{* \in Z}$,
\[\cdots F_{n+\Si (L)}^{(r)}HF^j(L, \p; P, J) \subset F_n^{(r)} 
HF^j(L, \p; P, J) \subset 
\cdots \subset I_j^{(r)}(L, \p; P, J),\] 
where $F_n^{(r)}HF^j(L, \p; P, J) = \ker (I_j^{(r)}(L, \p; P, J) \to
F_{n- \Si (L)}^{(r)}I_j^{(r)}(L, \p;P, J))$.
\end{lm}
\proof  The results follows from Definition~\ref{filter} and 
standard results
in \cite{sp} Chapter 9. \endproof

\begin{thm} \label{E1} For $\Si (L) \geq 3$,
there is a spectral sequence $(E^k_{n, j}, d^k)$ with
\[ E^1_{n,j}(L, \p; P, J) \cong I_n^{(r)}(L, \p; P, J), \ \ \ 
n\equiv j \ \ (mod \ \ \Si (L)), \]
\[ d^k\co E^k_{n, j}(L, \p; P, J) \to
E^k_{n + \Si (L) k +1, j +1}(L, \p; P, J) ,\]
and 
\[E^{\infty }_{n,j}(L, \p; P, J) \cong F_n^{(r)} HF^j(L, \p; P, J)/ 
F_{n + \Si (L)}^{(r)}HF^j(L, \p; P, J) . \]
In other words the spectral sequence $(E^k_{n, j}, d^k)$ converges to the 
${\BZ}_{\Si (L)}$-graded
symplectic Floer cohomology $HF^*(L, \p; P)$.
\end{thm}
\proof  Note that 
\[ F_n^{(r)} C_j(L, \p; P, J) /F_{n + \Si (L)}^{(r)} C_j(L, \p; P, J) = 
C_n ^{(r)}(L, \p; P, J) . \]
It is well--known from \cite{sp} that there exists a spectral sequence 
$(E^k_{n, j}, d^k)$ with $E^1$-- term given by the cohomology of $
F_n^{(r)} C_j(L, \p; P, J) /F_{n + \Si (L)}^{(r)} C_j(L, \p; P, J)$. So 
$E^1_{n,j}(L, \p; P, J) \cong I_n^{(r)}(L, \p; P, J)$
and $E^{\infty }_{n,j}(L, \p; P, J)$ is isomorphic to the bigraded 
${\BZ}_2$--module
associated to the filtration $F^{(r)}_*$ of the $\BZ$--graded symplectic
Floer cohomology $I_n^{(r)}(L, \p; P, J)$.
Note that the grading is unusual
(jumping by $\Si (L)$ in each step), we list the terms for $Z^k_{*, *}$
and $E^k_{*,*}$.
\[ Z^k_{n,j}(L, \p; P, J) = \{ x \in F_n^{(r)}C_j(L, \p; P, J) | \d x \in
F_{n+1 + \Si (L) k}^{(r)} C_{j+1}(L, \p; P, J) \}, \]
\[E^k_{n,j}(L, \p; P, J) = \]\[Z^k_{n,j}(L, \p; P, J)/
\{ Z^{k+1}_{n + \Si (L),j}(L, \p; P, J) +
\d Z^{k-1}_{n+ (k-1) \Si (L) -1,j-1}(L, \p; P, J)\}, \]
\[Z^{\infty}_{n,j}(L, \p; P, J) =
\{ x \in F_n^{(r)}C_j(L, \p; P, J)  | \d x = 0 \}, \]
\[E^{\infty}_{n,j}(L, \p; P, J) = \]\[Z^{\infty}_{n,j}(L, \p; P, J) /
\{ Z^{\infty}_{n + \Si (L),j}(L, \p; P, J) +
d Z^{\infty}_{n+ (k-1) \Si (L) -1,j-1}(L, \p; P, J)\}. \]
Thus 
$d^k\co E^k_{n, j}(L, \p; P, J) \to
E^k_{n + \Si (L) k +1, j +1}(L, \p; P, J)$ 
is induced from $\d$.

Since the Lagrangian intersections are transverse,
and $Z_{\p} = L \cap \p_1(L)$ is a finite
set, so the filtration $F^{(r)}_*$ is bounded
and complete from (\ref{coll}). Thus the spectral sequence converges to the
${\BZ}_{\Si (L)}$--graded symplectic Floer cohomology. 
\endproof

\noindent{\bf Remark}\qua For monotone Lagrangians $L$ with $\Si (L) \geq 3$
and $\p \in Symp_0(P, \o)$, the spectral sequence in Theorem~\ref{E1}
gives a precise information on $E^1_{*, *}$--term and all higher
differentials $d^k$ (see also later lemmae). The spectral sequence for
monotone Lagrangian $L$ with $\Si (L) \geq 2$ and $\p \in Ham (P, \o)$
$C^1$-sufficiently close to Id in Theorem IV of \cite{oh2}
gives a filtration from the Morse index of $L$. Our filtration is given by the 
integral lifting of the Maslov index. 

Theorem~\ref{E1} gives Theorem A (2). 
The following theorem is Theorem B.
\begin{thm} \label{44} 
For $\Si (L) \geq 3$, 
\begin{enumerate}
\item there exists an isomorphism
\[E^1_{n,j}(L, {\p}^0; P,J^0) \cong E^1_{n,j}(L, {\p}^1; P,J^1) ,\]
for any continuation $(J^{\lam }, {\p }^{\lam }) \in {\P }_{1, \ve /2}$
which is regular at ends.
\item for $k \geq 1$, $E^k_{n, j}(L, \p; P, J)$ is the symplectic invariant
under continuous deformations of $(J^{\lam }, {\p }^{\lam })$ within the
set of continuations.
\end{enumerate}
\end{thm}
\proof 
Clearly (2) follows from (1) by the Theorem 1 in \cite{sp} page 468.
By Theorem~\ref{E1}, there is an isomorphism:
\[E^1_{n,j}(L, {\p}^0; P,J^0) \cong I_n^{(r)}(L, {\p}^0; P,J^0).\]
By Theorem~\ref{invariant}, we have the isomorphism $$\p_{0,1}^n:
I_n^{(r)}(L, {\p}^0; P,J^0) \to I_n^{(r)}(L, {\p}^1; P,J^1),$$ 
which is compatible with the filtration. The isomorphism $\p_{0,1}^n$
induces an isomorphism on the $E^1$ term. So we obtain the result. \endproof

By Theorem~\ref{44}, $E^k_{n, j}(L, \p; P, J) = E^k_{n, j}(L, \p; P)$, 
$E^1_{n,j}(L,\p; P) = I_n^{(r)}(L,\p; P)$
are new symplectic invariants provided $\Si (L) \geq 3$.
All these new symplectic invariants should contain more information on
$(P, \o; L, \p)$. They are finer than the usual symplectic Floer cohomology 
$HF^*(L,\p; P) (* \in {\BZ }_{\Si (L)})$.
In particular, the minimal $k$ for which
$E^k(L, \p) = E^{\infty }(L, \p)$ should be meaningful, denoted by $k(L, \p)$.
The number $k(L, \p)$ is certainly a numerical invariant for the
monotone Lagrangian sub--manifold $L$ and $\p \in Symp_0(P)$. 

\begin{cor} \label{45}
For $\Si (L) \geq 3$ and $j \in {\BZ }_{\Si (L)}$, 
\[ \sum_{k \in Z} I^{(r)}_{ j + \Si (L) k }(L, \p; P) = HF^j (L, \p; P) \]
if and only if all the differentials $d^k$ in the spectral sequence 
$(E^k_{n, j}, d^k)$ are trivial, if and only if $k(L, \p) = 1$.
\end{cor} \endproof

In general, $\sum_{k \in Z} I^{(r)}_{ j + \Si (L) k }(L, \p; P)
\neq HF^j (L, \p; P)$ for $j \in {\BZ}_{\Si (L)}$. 
The $\BZ$--graded symplectic Floer cohomology $I_*^{(r)}(L,\p; P)$ 
can be thought as an integral lifting of the ${\BZ}_{\Si (L)}$--graded
symplectic Floer cohomology $HF^*(L, \p; P)$. From our construction of the 
$\BZ$--graded symplectic Floer cohomology, we have that

(i)\qua if $[r_0, r_1] \subset {\BR }_{L, \p}$, then
$E^{k, (r_0)}_{*, *}(L, \p; P) = E^{k, (r_1)}_{*, *}(L, \p; P)$;
 
(ii)\qua $E^{k, (r)}_{* + \Si (L), *}
(L, \p; P) = E^{k, (r+ \s (L))}_{* , *}(L, \p; P)$.

Since we work on field coefficients (over ${\BZ}_2$), we give a description
of $d^k$ in terms of generators.

\begin{lm} \label{mhd}
Modulo $E^{k-1}_{*, *}$, the differential $d^k\co E^k_{n, j} \to 
E^k_{n+k \Si (L)+1, j+1}$ is given by
\begin{equation} \label{hd}
d^k x = \sum_{- \mu^{(r)}(x) + \mu^{(r)}(y) = k \Si (L)+1}
\# {\hM}_J(x, y) \cdot y. 
\end{equation}
The differential $d^k$ extends linearly over $E^k_{n, j}$.
\end{lm}
\proof  
For $x \in E^k_{n, j}$, we have, by definition, that $x$ is survived
from all previous differentials. Since the coefficients are in a field
${\BZ}_2$, there are no torsion elements. So $x \in Z^k_{n, j}$,
$\mu^{(r)}(x) = n$ and $\d x \in F^{(r)}_{n+k \Si (L)+1}C_{j+1}(L, \p; P, J)$.
Thus 
\[d^kx = \d x = \sum_{y \in F^{(r)}_{n+k \Si (L)+1}C_{j+1}(L, \p; P, J)}
\# {\hM}_J(x, y) \cdot y. \]
The result follows. \endproof

\begin{lm} \label{col}
If $\mu^{(r)}_u(x, y) < (p+1)\Si (L)$ for any $u$ in 1--dimensional
moduli spaces of J--holomorphic curves with
$x, y \in C_*^{(r)}(L, \p; P, J)$, then $d^k = 0$ for $k \geq p+1$. So
the spectral sequence $E^*_{*,*}(L, \p; P, J)$ collapses
at least $(p+1)$th term.
\end{lm}
\proof 
Suppose the contrary. There is $x \in E^k_{n, j}$ such that
$d^k x \neq 0$ for $k \geq p+1$.
By Lemma~\ref{mhd}, 
\[\sum_{- \mu^{(r)}(x) + \mu^{(r)}(y) = k \Si (L)+1}
\# {\hM}_J(x, y) \cdot y \neq 0.\]
So there exist $y$ and $u\in {\hM}_J(x, y)$ such that there is
a nontrivial $J$--holomor\-phic curve $u$ with $\mu_u^{(r)} (x, y) = k \Si (L)+1$.
On the other hand, there exists a $J$--holomorphic curve $u$
with $\mu^{(r)}_u(x, y)\geq (p+1)\Si (L)+1$. Hence the
result follows from the contradiction. \endproof

By Lemma~\ref{col}, the lifted symplectic action actually measures the
how large the $J$--holomorphic curves by the integral lifted Maslov index.
If $a(x^{(r)}) - a(y^{(r)}) < (p+1) \s (L)$ for all $x, y 
\in C_*^{(r)}(L, \p; P, J)$, then $k(L, \p) \leq (p+1)$.

\begin{pro} For any compact monotone Lagrangian sub--manifold $L$
in $(P, \o)$ with
$\Si (L) \geq m +1, (m \geq 2)$, then
\begin{enumerate}
\item all the differentials $d^k$ are trivial for $k \geq 0$,
\item we have the following relation:
\[ \sum_{k \in Z}I^{(r_0)}_{j + \Si (L) k } (L, \p; P) = HF^j(L, \p; P) . \]
\end{enumerate}
\end{pro}
\proof  For the $\BZ$--graded symplectic Floer cochain complex
$C^{(r_0)}_*(L, \p_s;P,J)$, by Proposition~\ref{chek},
the Maslov index
satisfies the following:
\[0 < \max \mu^{(r_0)}(y) - \min \mu^{(r_0)}(x) =
\mu_{H_s}(y) - \mu_{H_s}(x) \leq m. \]
The result follows from the definition of $d^k$, $\Si (L) \geq m +1$
and Corollary~\ref{45}. \endproof

For any compact monotone Lagrangian embedding $L \subset {\BC}^m$,
we have $1 \leq \Si (L) \leq m$ and the
inequality is optimal based on Polterovich's
examples.

\begin{df} \label{plp}
(1)\qua The associated Poincar\'{e}--Laurent polynomial
$P(E^k, t)$ $(k \geq 1)$ of the spectral sequence is defined as:
\[P^{(r)}(E^k, t) = \sum_{n \in {\BZ}} (\dim_{{\BZ}_2} E^k_{n,j}) t^{n}. \]
(2)\qua The Euler number of the spectral sequence is defined to be the number
$\chi (E^k_{*,*}) = P^{(r)}(E^k, -1)$.
\end{df}

By Theorem~\ref{E1}, $P^{(r)}(E^1, t) = \sum_{n \in Z}
(\dim_{Z_2} I^{(r)}_n(L, \p ; P, J)) t^{n}$.
From Remark 5.b (ii), we have
\begin{equation} \label{shift}
P^{(r+\s (L))}(E^k, t) t^{\Si (L)} = P^{(r)}(E^k, t).
\end{equation}
We can compare two Poincar\'{e}--Laurent polynomials
$P^{(r)}(E^k, t)$ and $P^{(r)}(E^{\infty}, t)$ in the following.

\begin{pro} \label{pl}
For any monotone Lagrangian $L$ with $\Si (L) \geq 3$,
\[P^{(r)}(E^1, t) = \sum^k_{i =1} (1 + t^{-i \Si (L) -1}) \overline{Q}_i(t) +
P^{(r)}(HF^*, t), \]
where $k+1 = k(L, \p)$ and $\overline{Q}_i(t)$ ($i=1, 2, \cdots, k)$
are the Poincar\'{e}--Laurent polynomials of
nonnegative integer coefficients.
\end{pro}
\proof 
Let $Z^1_{n,j} = \ker \{d^1: E^1_{n,j} \to E^1_{n+ \Si (L) +1, j+1}\}$
and $B^1_{n,j} = \text{Im}\, d^1 \cap E^1_{n,j}$.
We have two short exact sequences:
\begin{equation} \label{zeb}
0 \to Z^1_{n,j} \to E^1_{n,j} \to B^1_{n+ \Si (L) +1, j+1} \to 0,
\end{equation}
\begin{equation} \label{bze}
0 \to B^1_{n,j} \to Z^1_{n,j} \to E^2_{n,j}  \to 0. 
\end{equation}
So the degree $\Si (L) +1$ of the differential $d^1$ derives the following.
\begin{equation} \label{esd}
P^{(r)}(E^1, t) = P^{(r)}(E^2, t) + (1+t^{- \Si (L) - 1}) P^{(r)}(B^1, t).
\end{equation}

Since the higher differential $d^i$ has degree $i \Si (L) +1$, we can repeat
(\ref{esd}) for $E^2$ and so on. Let $\overline{Q}_i(t) = P^{(r)}(B^i, t)$.
Note that $\op E^{\infty } \cong HF^*(L, \p; P)$
by Theorem~\ref{invariant} and Theorem~\ref{E1}. Thus
we obtain the desired result. \endproof

From the proof of Proposition~\ref{pl}, we have
\begin{equation} \label{jj}
P^{(r)}(E^l, t) = \sum^k_{i =l} (1 + t^{-i \Si (L) -1}) \overline{Q}_i(t) +
P^{(r)}(HF^*, t),
\end{equation}
for $1 \leq l \leq k$.
Since $\Si (L)$ is even, so we have
\begin{eqnarray*}
\chi (E^k_{*,*}) & = & \sum_{n \in Z} (-1)^n dim_{Z_2} E^k_{n,j} \\
  & = & \sum_{j =0}^{\Si (L) -1} (-1)^j dim_{Z_2} ({\bigoplus}_{n \equiv
j \pmod {\Si (L)}} E^k_{n,j}) .
\end{eqnarray*}
In particular, by Proposition~\ref{pl},
\begin{equation} \label{eule}
\chi (E^k_{*,*}) = \chi (E^{\infty}_{*,*}) = \chi (HF^*),
\ \ \ \ \mbox{for all $k \geq 1$}.
\end{equation}
For an oriented monotone Lagrangian sub--manifold $L_i$ in $(P_i, \o_i)$,
we have $I_{\o_i } = \lam I_{\mu, L_i}$ for the same $\lam \geq 0$. So the
product $L_1 \x L_2 \x \cdots \x L_s$ is also an
oriented Lagrangian in $(\prod_{i=1}^s P_i, \o_1 \op \o_2 \cdots \op \o_s)$.
For $u: (D^2, \bd D^2) \to (\prod_{i=1}^s P_i, \prod_{i=1}^s L_i)$, let
$p_i: \prod_{i=1}^s P_i \to P_i$ be the projection on the $i$-th factor
$(i=1, 2, \cdots, s)$. Now we obtain
\[I_{\mu \op \cdots \op \mu , L_1 \x  \cdots \x L_s}(u) = \sum_{i=1}^s
I_{\mu, L_i}(p_i u) , \]
which follows from the product symplectic form and the
K\"{u}nneth formula for the Maslov class in $\Lam ({\BC}^{m_1}) \x
\cdots \x \Lam ({\BC}^{m_s})$, where ${\BC}^{m_i} = (p_i\circ u)^*(T_xP_i)$
(see \cite {clm}). Hence we have
\begin{eqnarray*}
I_{\o_1 \op \cdots \op \o_s}(u) & = & \sum_{i=1}^s I_{\o_i } (p_i \circ u )  \\
 & = & \sum_{i=1}^s \lam I_{\mu, L_i}(p_i \circ u)  \\
 & = & \lam I_{\mu \op \cdots \op \mu , L_1 \x  \cdots \x L_s}(u) .
\end{eqnarray*}
So the product Lagrangian $L_1 \x \cdots \x L_s$ is also {\it monotone} in
$(\prod_{i=1}^s P_i, \o_1 \op \o_2 \cdots \op \o_s)$.
Note that $\la \Si (L_1 \x \cdots \x L_s) \ra = 
\la g.c.d (\Si (L_i): 1 \leq i \leq s) \ra$ (as an ideal in $\BZ$)
by the additivity of the
Maslov index. For simplicity, we will assume that $\Si (L_i) = \Si (L)
\geq 3$ for each $1 \leq i \leq s$.

Since $\Si (L_1 \x L_2) = \Si (L) \geq 3$ and $L_1 \x L_2$ 
is an oriented monotone
Lagrangian, we use the symplectic diffeomorphism 
$\p_1 \x \p_2 \in Symp_0(P_1 \x P_2)$. 
In particular, by Lemma~\ref{cs} and \cite{ch} we have
\[I^{(r_0)}_*(L_1 \x L_2, \p_1 \x \p_2; 
P_1 \x P_2) \cong H^{* + m_1+m_2 }(L_1 \x L_2 ;
{\BZ}_2) .\]
So there is a K\"{u}nneth formula for the
${\BZ}$--graded symplectic Floer cohomology.
\begin{equation} \label{kun}
I^{(r_0)}_*(L_1 \x L_2, \p_1 \x \p_2; P_1 \x P_2)
\cong I^{(r_0)}_*(L_1, \p_1; P_1) \ox I^{(r_0)}_*(L_2, \p_2; P_2).
\end{equation}
The torsion terms are irrelevant in this case due to the field 
${\BZ}_2$--coefficients. In terms of the filtration and
our spectral sequence, we have
\[E^1_{n,j}(L_1 \x L_2, \p_1 \x \p_2; P_1 \x P_2) \cong \]
\begin{equation} \label{eel}
\bigoplus_{n_1 + n_2 = n, j_1 + j_2 = j \pmod {\Si (L)} }
E^1_{n_1,j_1}(L_1, \p_1; P_1) \ox E^1_{n_2,j_2}(L_2, \p_2; P_2).
\end{equation}

\begin{pro} \label{ku}
For the monotone Lagrangian sub--manifold $L_1 \x L_2$,
there exists a spectral sequence
$E^k_{n,j}(L_1 \x L_2, \p_1 \x \p_2; P_1 \x P_2)$ which converges to
the ${\BZ}_{\Si (L)}$--graded symplectic Floer cohomology
$HF^*(L_1 \x L_2, \p_1 \x \p_2; P_1 \x P_2)$ with
\[E^1_{n,j}(L_1 \x L_2, \p_1 \x \p_2; P_1 \x P_2) =
I^{(r_0)}_*(L_1 \x L_2, \p_1 \x \p_2; P_1 \x P_2). \]
\end{pro}
\proof  The result follows from
the monotonicity of $L_1 \x L_2$, $\Si (L_1 \x L_2) = \Si (L) \geq 3$ 
and Theorem~\ref{E1}. \endproof

\begin{lm} \label{hd1}
For the spectral sequence $E^k_{n,j}(L_1 \x L_2, \p_1 \x \p_2; P_1 \x P_2)$
in Proposition~\ref{ku}, the higher
differential $d^k$ is given by
\[d^k|_{E^k_{n_1,j_1}(L_1, \p_1; P_1)} \ox 1 \pm  1 \ox
d^k|_{E^k_{n_2,j_2}(L_2, \p_2; P_2)} .\]
\end{lm}
\proof 
It is true for $k=0$ by (\ref{kun}). By the definition of $d^k$ (see
the proof of Theorem~\ref{E1}),
\[d^k \co  E^k_{n,j}(L_1 \x L_2, \p_1 \x \p_2; P_1 \x P_2) \to
E^k_{n+2Nk +1,j+1}(L_1 \x L_2, \p_1 \x \p_2; P_1 \x P_2),  \]
the term
$E^k_{n,j}$ is generated by cocycles in
$F_n^{(r,r)}C_j(L_1 \x L_2, \p_1 \x \p_2; P_1 \x P_2)$
(modulo $E^{k-1}_{n,j}$ term) since we use the field ${\BZ}_2$--coefficients.
After modulo $E^{k-1}_{*,*}$ we have the higher differential
\[d^k\co F^{(r)}_{n_1}C^{j_1}(L_1, \p_1; P_1) \ox
F^{(r)}_{n_2}C^{j_2}(L_2, \p_2; P_2)
\to \]\[\{F^{(r)}_{n_1+2Nk+1}C^{j_1+1}(L_1, \p_1; P_1) \ox
F^{(r)}_{n_2}C^{j_2}(L_2, \p_2; P_2)\} \]\[\oplus
\{F^{(r)}_{n_1}C^{j_1}(L_1, \p_1; P_1) \ox F^{(r)}_{n_2+2Nk+1}C^{j_2+1}
(L_2, \p_2; P_2) \} .\] Note that $d^k \circ d^k = 0$ by
the monotonicity and Proposition~\ref{ku}.
Thus the result follows
from the very definition of $d^k|_{E^k(L_i, \p_i; P_i)}$ for $i=1, 2$.
The sign is not important since the
coefficients are in ${\BZ}_2$.
\endproof

In general we can not expect the K\"{u}nneth formulae for the tensor product
of two spectral sequences due to the torsions (see \cite{li, sp}).
For our case with {\bf the field} ${\BZ}_2$--coefficients, we do have such a
K\"{u}nneth formula for the spectral sequence.
\begin{thm} \label{kus}
For the monotone Lagrangian $L_1 \x L_2$ in $(P_1 \x P_2, \o_1 \op \o_2)$ with
$I_{\o_i} = \lam I_{\mu, L_i}$ and $\Si (L_i) = \Si (L) \geq 3$ ($i=1, 2$), 
we have, for $k \geq 1$,
\[E^k_{n,j}(L_1 \x L_2, \p_1 \x \p_2; P_1 \x P_2) \cong \]\[
\bigoplus_{ n_1 + n_2 = n, j_1 + j_2 = j \pmod {\Si (L)}}
E^k_{n_1,j_1}(L_1, \p_1; P_1) \ox E^k_{n_2,j_2}(L_2, \p_2; P_2) .\]
\end{thm}
\proof 
The case of $k =1$ is the usual K\"{u}nneth formula for the cohomology
of $L_1 \x L_2$ (see (\ref{kun} and (\ref{eel}).
The result follows by induction on $k$ (see \cite{li, sp}) for
the field ${\BZ}_2$--coefficients.
\endproof

\begin{cor} \label{ksl}
{\rm(1)}\qua For a monotone Lagrangian $L \subset (P, \o)$, for $k \geq 1$ and
$s \geq 2$,
$E^k_{n,j}(L \x L \x \cdots \x L, \p \x \cdots \x \p; P^s) \cong$
\begin{equation*}
\bigoplus_{n_1 + \cdots + n_s = n, 
j_1 + \cdots + j_s = j \pmod {\Si (L)}}
\{\ox_{i=1}^s E^k_{n_i, j_i}(L, \p; P)\}. 
\end{equation*}
{\rm(2)}\qua The Poincar\'{e}--Laurent polynomial for $
E^k_{n,j}(L \x L \x \cdots \x L, \p \x \cdots \x \p; P^s)$ satisfies
\[P^{({\bf r})}(E^k(L \x L \x \cdots \x L, \p \x \cdots \x \p; P^s), t)
= \{P^{(r)}(E^k(L, \p; P), t)\}^s.\]
\end{cor}
\proof 
(1) follows from the induction proof of Theorem~\ref{kus} with 
$L_i = L (1\leq i \leq s)$ and the field
${\BZ}_2$--coefficients, and 
(2) follows from the definition of the Poincar\'{e}--Laurent
polynomial and (1). \endproof

\section{Applications}

\subsection{Hofer's energy and Chekanov's construction}

In this subsection, we relate our ${\BZ}$--graded symplectic Floer cohomology
with the one constructed by Chekanov 
\cite{ch}. We also show some results to illustrate the interactions among
the Hofer energy $e_H(L)$, the minimal symplectic action $\s (L)$ and the
minimal Maslov number $\Si (L)$.

Hofer \cite{ho1} introduced the notion of the disjunction energy or 
the displacement energy associated 
with a subset of symplectic manifold. 
Hofer's symplectic energy measures how large a variation of a (compactly 
supported) Hamiltonian function must be in order to push the subset off itself
by a time--one map of corresponding Hamiltonian flow. 
Hofer showed that the symplectic energy of every
open subset in the standard symplectic vector space is nontrivial.
See \cite{lm1} for more geometric study of the Hofer energy.

\begin{df} Let ${\cal H}$ be the space of compactly supported functions on
$[0,1] \x P$. The Hofer symplectic energy of a symplectic diffeomorphism 
$\p_1 \co P \to P$ is defined by
\begin{eqnarray*}
E(\p_1) &=& \inf \{ \int^1_0(\max_{x \in P} H(s, x) - \min_{x \in P} H(s, x))ds|
 \\
& & \mbox{$\p_1$ is a time one flow generated by $H \in {\cal H}$} \}. \\
e_H( L) &=& \inf \{ E(\p_1) : \p_1 \in Ham (P),  L \cap \p_1(L) = \emptyset 
\ \ \ \mbox{empty set} \}.\end{eqnarray*}
\end{df}

\begin{thm}\label{hse}{\rm\cite{ch}}\qua 
If $E(\p_1) < \s (L)$, $L$ is rational, and
$L$ intersects $\p_1(L)$ transversely, then
\[ \# (L \cap \p_1(L) ) \geq SB(L; Z_2), \]
where $SB(L; Z_2)$ is the sum of Betti numbers of $L$ 
with ${\BZ}_2$--coefficients.  
\end{thm}

\noindent{\bf Remark 5.a}\qua (1)\qua Polterovich 
used Gromov's figure $8$ 
trick and a refinement of Gromov's existence scheme of the $J$--holomorphic
disk to show that $e_H(L) \geq  \frac{\s (L)}{2}$. Chekanov extends the
result to $e_H(L) \geq \s (L)$ which is optimal for general Lagrangian
sub-manifolds (see \cite{ch} \S 1).  It is unclear whether 
Theorem~\ref{hse} remains true for the case of $E(\p_1) = \s (L)$. 

(2)\qua Sikorav showed that $e_H(T^m) \geq \s (T^m)$. The 
Theorem~\ref{hse}
generalizes Sikorav's result to all rational Lagrangian sub--manifolds.

Chekanov \cite{ch} used a restricted symplectic Floer cohomology in the study
of Hofer's symplectic energy of rational Lagrangian sub--manifolds. Denoted by
\begin{eqnarray*}
{\O }_s &= &\{ \g \in C^{\infty }([0,1], P) | 
\g (0) \in L, \g (1) \in \p_s (L)\}, \\
\O & =& \bigcup_{s \in [0,1]} {\O }_s  \subset [0,1] \x C^{\infty }([0,1], P). 
\end{eqnarray*}
Assume that $\p_s$ is generic so that $D=\{s\in [0, 1]: \p_s(L)\cap L 
\ \mbox{transversely}\}$ is dense in $[0,1]$.
One may choose the anti--derivative of $Da_s(z) \xi$ as $a_s\co {\O }_s \to 
{\BR }/{\s (L) {\BZ}}$, and fix $a_0$ with critical value $0$. 
Pick $z_s \in L \cap {\p}_s(L)$ such that modulo $\s (L) \BZ$, 
\begin{equation} \label{pf4}
0 < a_s (z_s) = \min \{a_s (x) \pmod {\s (L) {\BZ}} | x \in Z_{\p_s} \} 
< \s (L).\end{equation}
So it is possible to have $a_s (z_s)$ in
(\ref{pf4}) for $s\in D$ due to the constant factor for $a_s$.
Let $r_0 \neq 0$ be sufficiently small positive number in 
${\BR }_{L, {\p }_s} \cap (0, \s (L))$, say $0 < r_0 < \frac{1}{16} a_s(z_s)$. 
The condition $E(\p_1) < \s (L)$ provides that 
there is a unique $x \in Z_{{\p }_s}$ (the $x^{(r_0)}$) which
corresponds to the unique lift in $(r_0, r_0 + \s (L))$, 
\begin{equation} 0 < a_s(x) - a_s(z_s) < \s (L). 
\end{equation}
Under these restrictions, define the
free module $C_s$ over ${\BZ}_2$ generated by $L \cap \p_s (L)$ 
and the coboundary map
$\bd_s \in \text{End}(C_s)$ (see below) such that 
${\bd }_s \circ {\bd }_s =0$ \cite{ch}. Thus 
$H^*(C_s, {\bd }_s )$ is well--defined for every $s \in [0,1]$.
With the unique liftings of $x$ and $y$ in $(r_0, r_0+\s(L))$, we can identify 
Chekanov's restricted symplectic Floer cochain complex with our $\BZ$--graded
symplectic Floer cochain complex.

\begin{lm} \label{cs}
For $r_0$ as above, $C^{(r_0)}_*(L, \p_s; P, J) = C_s$.
Let ${\CM }_{J_s}(x, y)$ be the restricted moduli space of $J$--holomorphic 
curves $\{ u \in {\CM}(L, \p_s (L)) |
u^*(\o ) = a_s(x) - a_s(y) \}$. So 
\[ \bd^{(r_0)} x = \bd_s x = \sum \# {\hM }_{J_s}(x, y) y . \]
\end{lm} 
\proof 
Note that $\mu_u = \mu^{(r_0)}(y) - \mu^{(r_0)}(x) = dim {\CM}_{J_s}(x,y)$ 
for $u \in {\CM}_{J_s}(x,y)$. So $y$ in $\bd_s x$ is the element
in $C^{(r_0)}_{n+1}(L, \p_s; P,J)$; for any $u \in {\CM}_{J_s}(x,y)$ 
contributing in the coboundary of $\bd^{(r_0)} x$,
we have $u^*(\o ) = a_s(x) - a_s(y)$ by the 
Proposition 2.3 in \cite{fl2}. 
For the unique lifts in $(r_0, r_0 + \s (L))$ of $Z_{\p_s}$, 
the choice of $a_s$ gives arise to the 
one--to--one correspondence between ${\CM }_{J_s}(x, y)$ and
${\CM }_J(x, y)$ for $\mu^{(r_0)}(y) - \mu^{(r_0)}(x) = 1$. 
Therefore the coboundary maps agree on the ${\BZ}_2$--coefficients. \endproof

For $s$ sufficiently small, the point 
$x \in L \cap {\p }_s(L)$ is also a critical
point of the Hamiltonian function $H_s$ of $\p_s$. The Maslov index is 
related to the usual Morse index of the time--independent Hamiltonian $H$
with sufficiently small second derivatives:
\begin{equation}
\mu^{(r_0)}(x) = \mu_{H}(x) - m .
\end{equation}

\begin{pro} \label{chek}
For the $r_0$ as above, we assume that (i) $\Si (L) \geq 3$ and 
$E(\p_1) < \s (L)$, (ii) $L$ is monotone Lagrangian sub--manifold
in $P$, (iii) $L$ intersects
$\p_1(L)$ transversely. Then there is a natural isomorphism between 
\[I^{(r_0)}_*(L, \p_s; P) \cong H^{* + m}(L; Z_2)  \ \ \ \ 
\mbox{for $* \in Z$ and $s \in [0,1]$}. \]
\end{pro}
\proof 
For any $s, s^{'} \in [0,1]$ sufficiently close, 
by Theorem~\ref{invariant} we have  
\begin{equation} \label{ss}
I^{(r_0)}_*(L, \p_s; P) \cong I^{(r_0)}_*(L, {\p }_{s^{'}}; P) .
\end{equation}
For $s \in [0,1]$ sufficiently small, we have
\begin{equation} \label{so}
I^{(r_0)}_*(L, \p_s; P) \cong H^*(C_s, \bd_s) \cong H^{*+m}(L; Z_2).
\end{equation}
The first isomorphism is given by Lemma~\ref{cs} and the second by Lemma 3 of 
\cite{ch} and Theorem 2.1 of \cite{sz}. Then the result follows by applying
(\ref{ss}) and (\ref{so}). \endproof
\medskip

\noindent{\bf Remark 5.b}\qua (i)\qua 
Proposition~\ref{chek} provides the Arnold conjecture
for monotone Lagrangian sub--manifold with $\Si (L) \geq 3$ 
and $E(\p_1) < \s (L)$. Chekanov's result 
does not require the assumptions of the 
monotonicity of $L$ and $\Si (L) \geq 3$. 

(ii)\qua There is a natural relation 
$I^{(r_0 + \s (L))}_*(L, \p_s; P) = I^{(r_0)}_{* + \Si (L)}(L, \p_s; P)$ for
different choices of $a_s$ and $a_s + \s (L)$. In fact, 
for $s, s^{'} \in [0,1]$, 
\[I^{(r_0)}_*(L, \p_s; P) \cong I^{(r_0)}_*(L, {\p }_{s^{'}}; P). \]

For the Lagrangian sub--manifolds 
$L$ in $({\BC}^m, \o_0)$ with the standard symplectic
structure $\o_0 = - d \lam $. 
The 1--form $\lam $ is called {\em Liouville form},
which has the corresponding 
Liouville class $[\lam |_L] \in H^1(L, \BR)$. One of the fundamental results
in \cite{gr} is the non--triviality of the Liouville class.
 
\begin{thm}[Gromov \cite{gr}] \label{gromov} For any compact Lagrangian
embedding $L$ in ${\BC}^m$, 
the Liouville class $[\lam |_L] \neq 0 \in H^1(L, \BR)$.
\end{thm}

\begin{lm} \label{key}
If $E(\p_1) < (p+1) \s (L)$, 
then 
$k(L, \p) \leq (p+1)$ for the spectral sequence $E^k_{n, j}(L, \p; P, J)$.
\end{lm}
\proof  By definition of $E(\p )$, we have
\[a_{+, H} = \int_0^1 \max_{x \in P}H(s, x), \ \ \ \
a_{-, H} = \int_0^1 \min_{x \in P}H(s, x) , \]
for $\p_1$ a time--one flow generated by $H(s, x)$. For any nontrivial
$u \in {\hM}_J(x, y)$ which contributes in $d^k$, we have
\[I_{\o}(u) = a_u(x^{(r)}) - a_u(y^{(r)}) \leq a_{+, H} - a_{-, H} . \]
Thus for any $x, y \in C^{(r)}_*(L, \p; P, J)$, we have
$a_u(x^{(r)}) - a_u(y^{(r)}) \leq a_{+, H} - a_{-, H}$. So
\begin{equation}
\begin{split}
a_u(x^{(r)}) - a_u(y^{(r)}) & \leq \inf_{H} (a_{+, H} - a_{-, H}) \\
& = E(\p_1)  \\
& < (p+1) \s (L)\end{split} \end{equation}
By the monotonicity, $\lambda I_{\mu, \tilde{L}}(u) = I_{\o}(u) < (p+1)\s (L)$
($\mu^{(r)}(x, y) = I_{\mu, \tilde{L}}(u) < (p+1)\Si (L)$).
By Lemma~\ref{col}, we obtain the desired result. \endproof
 
\noindent{\bf Remark 5.c}\qua Note that $k(L, \p) = p+1$, one can not claim that
$E(\p) < (p+1)\s (L)$. For $k \geq p+1$, the property
$d^k = 0$ is the zero $\pmod 2$
of the 1--dimensional moduli space of $J$--holomorphic curves. There
are possible pairs $u_{\pm}$ of nontrivial $J$--holomorphic curves which have
$a_{u_{\pm}}(x^{(r)}) - a_{u_{\pm}}(y^{(r)}) \geq (p+1) \s (L)$.

There are four interesting numbers $\Si (L), \s (L), e_H(L)$ and $k(L, \p)$
of a monotone Lagrangian manifold $L$. 
Both of them intertwine and link
with the ${\BZ}$--graded symplectic Floer cohomology and its derived spectral
sequence. It would be interesting to study further relations among them.

\subsection{Internal cup--product structures}

For the monotone Lagrangian sub--manifold $L$ in 
$(P, \o)$ with $\Si (L) \geq 3$,
we have obtained an external product structure (cross product) in 
Theorem~\ref{kus}. In this subsection, we show that there is an internal
product structure on the spectral sequence and the symplectic Floer cohomology
of the Lagrangian sub--manifold $L$.

By a result of Chekanov \cite{ch} and Proposition~\ref{chek}, we can identify
$I^{(r)}_*(L, \p; P) \cong H^{*+m}(L; {\BZ}_2)$. Let $a$ be a cohomology class
in $H^p(L; {\BZ}_2)$. Define a map
\begin{equation} \label{msi}
a \cup \co C^{(r)}_n(L, \p; P, J) \to C^{(r)}_{n+p}(L, \p; P, J)
\end{equation}
\[ x \mapsto \sum_{- \mu^{(r)}(x) + \mu^{(r)}(y) = p}
\# ({\CM}_J(x, y) \cap i_*(PD_L(a))) \cdot y, \]
where $i\co 
L \hookrightarrow P$
is the Lagrangian imbedding, $\# ({\CM}_J(x, y) \cap i_*(PD_L(a)))$ is the 
algebraic number of $i_*(PD_L(a))$ intersecting ${\CM}_J(x, y)$,
and $PD_L(a)$ is the Poincar\'{e} dual of $a$ in $L$.

\begin{pro} \label{wedc}
The map $a \cup $ in (\ref{msi}) is well--defined for the monotone
Lagrangian embedding, and $\bd^{(r)}_{n+p} \circ (a \cup ) = (a \cup ) \circ
\bd^{(r)}_n$.
\end{pro}
\proof  Note that $i_*(PD_L(a))$ is a divisor in $P$ where the 
intersection ${\CM}_J(x, y)$ $\cap\, i_*(PD_L(a))$ can be made transversally without
holomorphic bubblings (see \cite{fl2, fl3} Theorem 6). In fact, we can apply
the similar argument in Proposition 4.1 of \cite{li1} to the ${\BZ}$--graded
symplectic Floer cohomology of the Lagrangian $L$ in order to avoid the bubbling
issue. Now it suffices to check 
the map $a \cup$ commutes with the ${\BZ}$--graded differential.

Without holomorphic bubblings, the partial compactification of ${\CM}_J(x, z)$
with only $k$--tuple holomorphic curves can
be described as
\[\ov{{\CM}_J(x, z)} =
\cup (\x_{i=0}^{k-1} {\CM}_J(c_{i}, c_{i+1})),\]
the union over all sequence $x=c_0, c_1, \cdots,
c_k=z$ such that
${\CM}_J(c_{i}, c_{i+1})$ is nonempty for all
$0 \leq i \leq k-1$.
For any sequence $c_0, c_1, \cdots, c_k \in \text{Fix}\, (\p)$,
there is a gluing map
\[G\co \x_{i=0}^{k-1} {\hM}_J(c_i, c_{i+1}) \x
\D^k \to \ov{{\CM}_J(x,z)},\]
where $\D^k=\{(\lam_1, \cdots, \lam_k) \in [-\infty, \infty]^k:
1+\lam_{i} < \lam_{i+1}, 1 \leq i \leq k -1 \}$.
\begin{enumerate}
\item The image of $G$ is a neighborhood of
$\x_{i=0}^{k-1} {\hM}_J(c_i, c_{i+1})$ in the
compactification with only $k$--tuple holomorphic curves.
\item The restriction of $G$ to
$\x_{i=0}^{k-1} {\hM}_J(c_i, c_{i+1}) \x \text{Int}\,
(\D^{k})$ is a diffeomorphism onto its image.
\item The extension of the
gluing map is independent of $u_1, \cdots, u_{j-1}$ and $u_{j+p+1}, \cdots,
u_k$ provided $\lam_{j-1} = - \infty$ and $\lam_{j+p+1} = + \infty$.
$$G(u_1, \cdots, u_k, - \infty, \cdots, - \infty, \lam_j, \cdots,
\lam_{j+p}, +\infty, \cdots, + \infty) $$
$$=
G(u_j, \cdots, u_{j+p}, \lam_j, \cdots, \lam_{j+p}).$$
\end{enumerate}
For $x \in C_n^{(r)}(L, \p; P, J)$ and
$z \in C_{n+p+1}^{(r)}(L, \p; P, J)$, the space
$K = {\CM}_J(x, z)\cap i_*(PD(a))$
is a 1--dimensional manifold in $\O_{\p}$ with $p = \deg (a)$.
We have
\begin{equation*}
\begin{split}
0  & = \int_{\ov{{\CM}_J(x,z)} } d PD_P^{-1}(i_*(PD_L(a)))
\\ & = \int_{\bd \ov{{\CM}_J(x,z)} } PD_P^{-1}(i_*(PD_L(a)))\\ 
& = 
\bd \ov{{\CM}_J(x,z)} \cap i_*(PD_L(a)),
\end{split} \end{equation*}
where $PD_X$ stands for the Poincar\'{e} dual of the space $X$.
Since $G$ is a local diffeomorphism near
$\bd \ov{{\CM}_J(x,z)}$, we can integrate
$PD_P^{-1}(i_*(PD_L(a)))$ (by changing variables) over 
$\x_{i=0}^{k-1}{\hM}_J(c_{i}, c_{i+1}) \x \D^k$. We have
$$\bd (\x_{i=0}^{k-1}{\hM}_J(c_{i},c_{i+1}) \x \D^k) =
\x_{i=0}^{k-1}{\hM}_J(c_{i}, c_{i+1}) \x \bd \D^k.$$
From the definition of $\D^k$, we get
$\bd (\D^k) = \D_+^{k-1} \coprod - \D_-^{k-1}$,
where $\D_-^{k-1} = \{(- \infty, \lam_2, \cdots, \lam_k) \in \D^{k}\}$ and
$\D_+^{k-1} = \{(\lam_1, \cdots, \lam_{k-1}, + \infty) \in \D^{k}\}$. Thus
\begin{equation} \label{stok1}
\begin{split}
0 & = \int_{\bd \ov{{\CM}_J(x,z)} } PD_P^{-1}(i_*(PD_L(a))) \\
& = \la \bd G^{-1}(\x_{i=0}^{k-1}{\hM}_J(c_{i},
c_{i+1}) \x \D^{k}), PD_P^{-1}(i_*(PD_L(a))) \ra \\
& = \la \x_{i=0}^{k-1}{\hM}_J(c_{i},
c_{i+1}) \x \bd (\D^{k}), G^*(PD_P^{-1}(i_*(PD_L(a))))\ra \\
& = \la \x_{i=0}^{k-1}{\hM}_J(c_{i}, c_{i+1}) \x \D_+^{k-1}, 
G^*(PD_P^{-1}(i_*(PD_L(a))))\ra \\
&
- \la \x_{i=0}^{k-1}{\hM}_J(c_{i}, c_{i+1}) \x \D_-^{k-1}, 
G^*(PD_P^{-1}(i_*(PD_L(a))))\ra.\\
\end{split} \end{equation}
The image of $G$ is independent of
$u_1$ in $\D_-^{k-1}$ and $u_k$ in $\D_+^{k-1}$ respectively.
For $\lam_1 \to - \infty$, $$G(\x_{i=1}^{k-1}{\hM}_J(c_i, c_{i+1})
\x \D_-^{k-1}) = \ov{{\CM}_J(c_1, z)}.$$
We have the dimension counting as follows from
$\sum_{i=0}^{k-1} (\mu^{(r)}(c_{i}) - \mu^{(r)}(c_{i+1})) = p+1$,
\[\dim \ov{{\CM}_J(c_1, z)} =
(p+1) - (\mu^{(r)}(x) - \mu^{(r)}(c_1)).\]
By the transversal of the intersection with $i_*(PD_L(a))$ in $\O_{\p}$,
the only possible nontrivial contribution of
$\ov{{\CM}_J(c_1, z)} \cap i_*(PD_L(a))$ is from
$\dim \ov{{\CM}_J(c_1, z)} = p$. Hence
$(p+1) - (\mu^{(r)}(x) - \mu^{(r)}(c_1))=p$ if and only if
$\mu^{(r)}(x) - \mu^{(r)}(c_1) = 1$. Therefore we obtain
\[\la \x_{i=0}^{k-1}{\hM}_J(c_{i},
c_{i+1}) \x \D_-^{k-1}, G^*(PD_P^{-1}(i_*(PD_L(a))))\ra \]\[=
\# {\hM}^1_J(x, c_1) \cdot
\# {\CM}_J(c_1, z) \cap i_*(PD_L(a)),\] this gives the term
$(a\cup ) \circ \bd_n^{(r)}(x)$. Similarly for 
$\lam_k \to + \infty$,
\begin{equation} \label{stok3}
\la \x_{i=0}^{k-1}{\hM}_J(c_{i},
c_{i+1}) \x \D_+^{k-1}, G^*(PD_P^{-1}(i_*(PD_L(a))))\ra = 
\bd_{n+p}^{(r)} \circ (a \cup )(x).
\end{equation}
Hence the result follows. \endproof

Now the map $a \cup$ defined in (\ref{msi}) induces a map
(still denoted by $a \cup$) on the $\BZ$--graded symplectic Floer cohomology by
Proposition~\ref{wedc},
\begin{equation} \label{fim}
a \cup \co I^{(r)}_n (L, \p; P, J) \to I^{(r)}_{n+p}(L, \p; P, J) .
\end{equation}
Since $H^*(L; {\BZ}_2)$ is a graded algebra with cup
product as multiplication, the $H^*(L; {\BZ}_2)$--module structure of
$I^{(r)}_* (L, \p; P, J)$ is given by the following commutative diagram:
\[\begin{array}{ccc}
H^*(L; {\BZ}_2) \ox H^*(L; {\BZ}_2) \ox I^{(r)}_* (L, \p; P, J) &
\stackrel{1 \ox  \psi}{\longrightarrow}&
H^*(L; {\BZ}_2) \ox I^{(r)}_* (L, \p; P, J) \\
\Big\downarrow\vcenter{%
\rlap{$\cup \ox 1$}} & &
\Big\downarrow\vcenter{%
\rlap{$\psi $}}\\
H^*(L; {\BZ}_2) \ox I^{(r)}_* (L, \p; P, J) &
\stackrel{\psi}{\longrightarrow} &  I^{(r)}_* (L, \p; P, J)
\end{array} \]
where $\psi \co H^*(L; {\BZ}_2) \ox I^{(r)}_* (L, \p; P, J) \to 
I^{(r)}_* (L, \p; P, J)$
is given by $\psi (a , \cdot ) = a \cup \cdot$ in (\ref{fim}). Thus
we have $\{(a \cup ) \circ (b \cup)\}(x)$ equals
\[\sum_{y\in C^{(r)}_{n+\deg (b)}} \# {\CM}_J(x,y) \cap i_*(PD_L(b)) \cdot
\sum_{z\in C^{(r)}_{n+\deg (b)+\deg (a)}} \# {\CM}_J(y, z) \cap
i_*(PD_L(a)) \cdot z.\]

\begin{pro} \label{ims}
If $i\co L \hookrightarrow P$ 
induces a surjective map $i^*\co H^*(P; {\BZ}_2) \to
H^*(L; {\BZ}_2)$, then the $\BZ$--graded symplectic Floer cohomology 
$I^{(r)}_*(\p, P)$ has an
$H^*(L; {\BZ}_2)$--module structure.
\end{pro}
\proof  If $i^*\co H^*(P; {\BZ}_2) \to
H^*(L; {\BZ}_2)$ is surjective, then $i_* \circ PD_L \circ i^* = PD_P$.
Thus $i_*(PD_L(a \cup b)) = i^*(PD_L(i^*(a_P \cup b_P))) = PD_P(a_P \cup
b_P)$ for some classes $a_P, b_P \in H^*(P; {\BZ}_2)$. Therefore
\begin{equation}
\begin{split}
\# {\CM}_J(x, z) \cap i_*(PD_L(a \cup b)) & = \la PD_P^{-1}(i_*(PD_L(a \cup b)))
, {\CM}_J(x, z) \ra \\  & = \la a_P \cup
b_P, {\CM}_J(x, z) \ra \\
& = \la D^*(a_P \x b_P), {\CM}_J(x, z) \ra \\
& = \la a_P \x b_P  , \D_* {\CM}_J(x, z) \ra, \\
\end{split} \end{equation}
where $\D_*$ is the unique chain map up to chain homotopy induced from
the diagonal map $\D\co \O_{\p} \to \O_{\p} \x \O_{\p}$, 
and $H^*(P; {\BZ}_2)$ is
viewed as a subring of $H^*({\O}_{\p}; {\BZ}_2)$ 
(see \cite{fl3} \S 1c or \cite{li1}).
By the monotonicity and the energy formula, for any elements in
${\CM}_J(x,y) \x {\CM}_J(y, z)$, we have
\[ r < \ti{a}_J(x^{(r)}) \leq \ti{a}_J(y^{(r)}); \ \ \ \ \
\ti{a}_J(y^{(r)}) \leq \ti{a}_J(z^{(r)}) < r + 2\a N. \]
So by the uniqueness of the lifting $y^{(r)}$,
we have $y \in C^{(r)}_*(L, \p; P, J)$.
For $\mu^{(r)}(y) \neq n + deg (b)$, then, by the dimension counting,
\[ \la a_P \x b_P  , {\CM}_J(x,y) \x {\CM}_J(y,z) \ra = 0.\]
Thus we obtain $\la a_P \x b_P  , \cup_{y\in C^{(r)}_{n+\deg (b)}}
{\CM}_J(x,y) \x {\CM}_J(y,z) \ra  $
\begin{equation} \begin{split}
& =
\sum_{\mu^{(r)}(y) =n + deg (b)} \la b_P, {\CM}_J(x,y) \ra \cdot \la
a_P, {\CM}_J(y,z) \ra \\
& = \sum_{\mu^{(r)}(y) =n + deg (b)} \la 
PD_P^{-1}(i_*PD_L(b)), {\CM}_J(x,y) \ra \cdot \la
PD_P^{-1}(i_*PD_L(a)), {\CM}_J(y,z) \ra \\
& = \sum_{y\in C^{(r)}_{n+\deg (b)}}\# ({\CM}_J(x, y) \cap i_*(PD_L(b)))
\cdot \# ({\CM}_J(y, z) \cap i_*(PD_L(a))) z.\\
\end{split} \end{equation}
This equals to $\la (a \cup ) \circ (b \cup ) (x), z\ra$.
\endproof

One can verify that the module structure is invariant under the
compact continuation $(J^{\lam}, \p^{\lam}) \in {\P}_{1, \ve/2}$
by the standard method in \cite{fl3, li1}. Note that one can define the
action of $H^*(L; {\BZ}_2)$ for any monotone Lagrangian sub--manifold $L$,
but the module structure requires for the special property of the Lagrangian
embedding $i\co L \hookrightarrow P$ in Proposition~\ref{ims}.

Hence we obtain a $H^*(L; {\BZ}_2)$--module structure on $I^{(r)}_*(L, \p; P)$.
Now we associate the internal product structure on the ${\BZ}$--graded
symplectic Floer cohomology:
\begin{equation} \label{piz}
I^{(r)}_{n_1-m}(L, \p; P) \x I^{(r)}_{n_2-m}(L, \p; P) \stackrel{\cup}{\to}
I^{(r)}_{n_1+n_2-m}(L, \p; P),
\end{equation}
as a bilinear form defined by 
$I^{(r)}_{n_1-m}(L, \p; P) \cong H^{n_1}(L; {\BZ}_2)$ and (\ref{fim}).
Note that the index--shifting makes the compatibility of the usual cup--product
of the cohomology ring on $H^*(L; {\BZ}_2)$: $H^{n_1}(L; {\BZ}_2) \x
H^{n_2}(L; {\BZ}_2) \stackrel{\cup}{\to}
H^{n_1+n_2}(L; {\BZ}_2)$.

\begin{cor} \label{inh} If $i\co L \hookrightarrow P$ induces
a surjective map $i^*\co H^*(P; {\BZ}_2) \to H^*(L; {\BZ}_2)$ and
$\Si (L) \geq 3$, then
\[H^*(L; {\BZ}_2) \to End(I^{(r)}_*(L, \p, P; {\BZ}_2))\]
is an injective homomorphism.
\end{cor}

Corollary~\ref{inh} generalizes Theorem 3 of Floer \cite{fl1} of 
$\pi_2(P, L) = 0$ to the case of monotone Lagrangian sub-manifolds $L$ with
$\Si (L) \geq 3$ if the intersection $L$ with $\p (L)$ is transverse,
where $\p \in Symp_0(P)$ not necessary $\p \in Ham(P)$.
Note that one may combine our construction in \S 2 - \S 4 with the one in 
\S 3 of \cite{fl1} for the non-transverse points in $L \cap \p (L)$.
With the gluing result of trajectories along those degenerate points in 
$L \cap \p (L)$, we obtain that for any exact Hamiltonian $\p$ on the
monotone pair $(P, L; \o)$ with $\Si (L) \geq 3$, the number
$\# L \cap \p (L)$ of intersections of $\p (L)$ with $L$ is greater than or
equal to the ${\BZ}_2$-cuplength of $L$ via (\ref{piz}) and Corollary~\ref{inh}.
I.e., a generalization of Theorem 1 of \cite{fl1} is achieved for a 
symplectic manifold $(P, \o)$ and its monotone Lagrangian sub-manifold
$L$ without the hypothesis $\pi_2 (P, L)  = 0$.

\begin{lm} \label{fpre}
For the Lagrangian embedding in Proposition~\ref{ims}, the map
$a \cup $ defined in (\ref{msi}) 
preserves the filtration in Definition~\ref{filter}, and $d^k \circ (a \cup)
= (a \cup ) \circ d^k$, where $d^k$ is the higher differential  in the spectral
sequence of Theorem~\ref{E1}.
\end{lm}
\proof  Recall $F_n^{(r)}C^j(L, \p, P; J) 
= \sum_{k \geq 0}C^{(r)}_{n+ \Si (L)k}(L, \p, P;J)$.
Thus by our definition $a \cup$ in (\ref{msi}), we have
\[a \cup : \sum_{k \geq 0}C^{(r)}_{n+\Si (L)k}(L, \p, P;J) \to
\sum_{k \geq 0}C^{(r)}_{n+\Si (L)k+\deg (a)}(L, \p, P;J),  \]
which induces a map $a \cup \co F_n^{(r)}C^j(L, \p, P; J) \to
F_{n+\deg (a)}^{(r)}C^{j_a}(L, \p, P; J)$ with 
$j_a \equiv j+ \deg (a) \pmod {\Si (L)}$.
Thus we obtain a filtration preserving homomorphism with degree $\deg (a)$.
Element in $E^k_{n,j}$ is a survivor from previous differentials, and
$x \in E^k_{n,j}$ is an element in $F_n^{(r)}C^j(L, \p, P, J)$ with $\d x
\in F^{(r)}_{n+1+\Si (L)k}C^{j+1}(L, \p, P, J)$. The differential $d^k$ (induced
from $\d$) is counting the signed one--dimensional moduli space of
$J$--holomorphic curves from $x \in E^k_{n,j}$ to 
$y \in E^k_{n+1+\Si (L)k, j+1}$.
The diagram
\begin{equation} \label{edk}
\begin{array}{ccc}
E^k_{n,j} & \stackrel{d^k}{\longrightarrow} & E^k_{n+1+\Si (L)k, j+1}\\
\downarrow {a \cup} & & \downarrow {a \cup} \\
E^k_{n+\deg (a), j_a} & \stackrel{d^k}{\longrightarrow} &
E^k_{n+1+\Si (L)k+\deg (a), j_a+1}
\end{array} \end{equation}
is commutative by the same method of the proof in Proposition~\ref{wedc}.
\endproof

\begin{thm} \label{cupp}
For the monotone Lagrangian $L \stackrel{i}{\hookrightarrow} P$
with $i^*\co H^*(P; {\BZ}_2) \to H^*(L; {\BZ}_2)$ 
surjective and $\Si (L) \geq 3$,
the spectral sequence $E^k_{*-m, *}(L, \p; P)$ has an ring structure
which is descended from the cohomology ring $(H^*(L; {\BZ}_2), \cup )$.
\end{thm}
\proof  As we see from the above, there is an internal product 
structure on $E^1_{*-m, *}$: $H^{n_1}(L; {\BZ}_2) \x E^1_{n_2-m, j_2}(L, \p; P) =$
\begin{equation} \label{e1ps}
E^1_{n_1-m, j_1}(L, \p; P) \x E^1_{n_2-m, j_2}(L, \p; P) \to 
E^1_{n_1+n_2-m, j}(L, \p; P),
\end{equation}
where $E^1_{n_1-m, j_1}(L, \p; P) \cong H^{n_1}(L; {\BZ}_2)$ and
$j \equiv j_1+j_2+m \pmod {\Si (L)}$. Hence the module structure in 
Proposition~\ref{ims} shows that the internal product is associative through
the identifications. Furthermore the module structure descends to the spectral
sequence by Lemma~\ref{fpre}:
$$H^{n_1}(L; {\BZ}_2) \x E^k_{n_2-m, j_2}(L, \p; P) \stackrel{\cup }{\to}
E^k_{n_1+n_2-m, j}(L, \p; P).$$
As a subgroup $E^k_{n_1-m, j_1}(L, \p; P) \subset H^{n_1}(L; {\BZ}_2)$, 
elements in $E^k_{n_1-m, j_1}(L, \p; P)$ are
the survivors which are neither an image of some $d^k$ nor the source
of the nontrivial $d^k$ due to the ${\BZ}_2$--coefficients. Hence they 
are exactly the cohomology classes
in $H^{n_1}(L; {\BZ}_2)$.
In general, the $k$-th term $E^k_{n_1-m, j_1}(L, \p; P)$ lies in the quotient of
$E^{k-1}_{n_1-m, j_1}(L, \p; P)$ which is not necessary a subgroup of
$H^{n_1}(L; {\BZ}_2)$. Therefore
there is an associate (from the module structure)
internal product structure on $E^k_{*-m, *}(L, \p; P)$:
\begin{equation}\label{uuu}
 E^k_{n_1-m, j_1}(L, \p; P) \x E^k_{n_2-m, j_2}(L, \p; P) 
\stackrel{\cup_k}{\to} E^k_{n_1+n_2-m, j}(L, \p; P),\end{equation}
via the identification of the first factor
$E^k_{n_1-m, j_1}(L, \p; P) \subset H^{n_1}(L; {\BZ}_2)$.
\endproof

In particular, $(E^{\infty}_{*-m, *}(L, \p; P), \cup_{\infty}) =
(HF^{*-m}(L, \p; P), \cup_{\infty})$ is an associative ring on the
${\BZ}_{\Si (L)}$--graded symplectic Floer cohomology which is induced
from the cohomology ring of the imbedded monotone Lagrangian
sub--manifold $L$. Note that the index shifting (\ref{uuu}) in the internal
product is not the usual one due to the relation between the Maslov index and
the Morse index.

For $\p \in \text{Symp}_0(P)$ and the monotone Lagrangian embedding
$i\co L \hookrightarrow P$ with $i^*\co H^*(P; {\BZ}_2) \to H^*(L; {\BZ}_2)$
surjective and $\Si (L) \geq 3$, we obtain the ring structure on 
$E^k_{*-m, *}(L, \p; P)$ for every $k\geq 1$.
From the quantum effect arisen from the higher differentials in the
spectral sequence in Theorem~\ref{E1}, the ring
$(HF^{*-m}(L, \p; P), \cup_{\infty})$ can be thought of as
the quantum effect of the regular cohomology ring
$(H^*(L; {\BZ}_2), \cup)= (E^1_{*-m, *}(L, \p; P), \cup_1)$
of the monotone Lagrangian sub--manifold $L$ with embedding
$i\co L \to P$.

\noindent{\bf Remark}\qua There is a general approach for the
$A_{\infty}$--structure on the symplectic Floer cohomology of Lagrangians
with Novikov ring coefficients in \cite{fooo}. The $A_m$--multiplicative 
structure is defined 
by pair--of--pants construction which is quite a complicated
and hard
in terms of computation. Our cup--product is induced from the
usual cup-product of $H^*(L; {\BZ}_2)$ and incorporated with those
higher differentials. Our construction and the cup-product structure
have more algebraic topology techniques in terms of computation.

\noindent{\bf Acknowledgments}\qua The author is very grateful to the referee
for his/her comments and corrections on a previous version of the paper.
The research of the author is partially supported by NSF grant number
DMS-0245323.

\Addresses\recd
\end{document}